\documentclass
[12pt]{amsart}
\usepackage{amssymb}
\usepackage[all]{xy}

\newcommand{\Q}{{\mathbb {Q}}}

\newcommand{\R}{{\mathbb{R}}}
\newcommand{\Z}{{\mathbb{Z}}}
\newcommand{\C}{{\mathbb{C}}}

\newcommand{\Ad}{{\operatorname{Ad}}}

\newcommand{\SL}{\operatorname{SL}}

\newcommand{\OO}{{\mathcal O}}

\newcommand{\rank}{{\rm rank}}

\theoremstyle{plain}
\newtheorem{thm}{Theorem}[section]
\newtheorem{lem}[thm]{Lemma}
\newtheorem{prop}[thm]{Proposition}
\newtheorem{cor}[thm]{Corollary}

\theoremstyle{definition}

\def\bydefn{\stackrel{def}{=}}

\title[Locally divergent orbits]{Locally divergent orbits on Hilbert modular spaces}

\author{George Tomanov}
\address{Institut Camille Jordan, Universit\'e Claude Bernard - Lyon
I, B\^atiment de Math\'ematiques, 43, Bld. du 11 Novembre 1918,
69622 Villeurbanne Cedex, France {\tt tomanov@math.univ-lyon1.fr}}

\begin{document}

\begin{abstract}
We describe the closures of locally divergent orbits under the
action of tori on Hilbert modular spaces of rank $r \geq 2$. In
particular, we prove that if $D$ is a maximal $\R$-split torus
acting on a real Hilbert modular space then every locally divergent
non-closed orbit is dense for $r>2$ and its closure is a finite
union of tori orbits for $r=2$. Our results confirm an orbit
rigidity conjecture of G.A.Margulis in all cases except for (i) $r = 2$
and, (ii) $r > 2$ and the Hilbert modular space corresponds to a
$\mathrm{CM}$-field; in the cases (i) and (ii) our results
contradict the conjecture.

As an application, we describe the set of values at integral points
of collections of non-proportional, split, binary, quadratic forms
over number fields.

\end{abstract}

\maketitle

\section{Introduction}

During the last decade the problems of the descriptions of orbit
closures for actions of maximal split tori on homogeneous spaces
appear to be among the central ones in homogeneous dynamics. The
interest in such problems is motivated to a large extent by number theoretic
applications. One example about the efficiency of the homogeneous dynamics approach in
the number theory is G.A.Margulis' proof of the long-standing Oppenheim conjecture
 \cite{Margulis Oppenheim 2}. In our days this
approach looks quite promising regarding the still open Littlewood
conjecture \cite[\S2]{Margulis Oppenheim}. We refer to \cite{survey}
and \cite{Margulis problems survey} for an account of results and
conjectures on the subject. In the present
paper{\footnote{\textit{To appear in "International Mathematics
Research Notices", 2012.}}}, as an application of the main results,
we give an explicit description of the set of values at integral
points of a collection of non-proportional, split, binary quadratic
forms over number fields.

Let us introduce some basic terminology. Let $K$ be a number field,
$\OO$ its ring of integers and $K_i, \ 1 \leq i \leq r,$ all the
archimedean completions of $K$. {\it Throughout this paper we assume
that $r \geq 2$.} Put $G = \underset{i=1}{\overset{r}{\prod}}G_i$,
where $G_i = \SL(2,K_i)$, and let $\Gamma = \SL(2,\OO)$ be
identified with its image in $G$ under the diagonal embedding. The
subgroup $\Gamma$ is a non-uniform irreducible lattice in $G$ and by
the arithmeticity theorem (cf. \cite{Margulis arithmeticity 1},
\cite{Margulis arithmeticity 2}, \cite{Selberg}) all non-uniform
irreducible lattices in $G$ arise by this construction up to
conjugation and commensurability. The quotient space $G/\Gamma$ is
called {\it Hilbert modular space of rank $r$}. Denote by $\pi:G
\rightarrow G/\Gamma$ the natural projection. Let $D_i$ be the
connected component of the identity of the diagonal subgroup of
$G_i$ and let $D_{i,\R}$ be the connected component of the identity
of the subgroup of {\it real} matrices in $D_i$. So, $D_{i,\R} =
D_i$ if $K_i = \R$. For every non-empty $I \subset \{1, \cdots ,r\}$
we denote $D_I = \underset{i \in I}{\prod}D_i$ and $D_{I,\R} =
\underset{i \in I}{\prod}D_{i,\R}$. When $I = {\{1, \cdots ,r\}}$ we
write $D$ and $D_{\R}$ instead of $D_I$ and $D_{I,\R}$,
respectively. By a torus (respectively, an $\R$-split torus or,
simply, a split torus) in $G$ we mean a subgroup conjugated to a
closed connected subgroup of $D$ (respectively, $D_{\R}$). An orbit
$D_I\pi(g)$ is called {\it locally divergent} if $D_i\pi(g)$ is
divergent for all $i \in I$. (Recall that if $H$ is a closed
non-compact subgroup of $G$ and $x \in G/\Gamma$ then the orbit $Hx$
is divergent if the orbit map $h \mapsto hx$ is proper or,
equivalently, if $\{h_nx\}$ leaves compact subsets of $G/\Gamma$
whenever $h_n$ leaves compact subsets of $H$.) The orbit
$D_{I,\R}\pi(g)$ is locally divergent if and only if the orbit
$D_I\pi(g)$ is locally divergent. The description of the divergent
$D_i$-orbits (and, therefore, the divergent $D_{i,\R}$-orbits)
follows from the general results of \cite{Toma} (see \S\ref{locally
divergent orbits}).

The following conjecture is a special case of
a conjecture of G.A.Margulis \cite[Conjecture 1]{Margulis problems survey}.
\medskip

\textbf{Conjecture A}(\emph{orbit rigidity}): If $\#I \geq 2$ then
every orbit $D_{I,\R}x$, $x \in G/\Gamma$, has {\it homogeneous
closure}, that is, $\overline{D_{I,\R}x} = Fx$, where $F$ is a
closed subgroup in $G$ containing $D_{I,\R}$.

\medskip

Broadly speaking, the general \cite[Conjecture 1]{Margulis problems survey}
says that the closure of an orbit for the action of an $\R$-split
torus $T$ of dimension $\geq 2$ on a homogeneous space of finite
volume $G/\Gamma$ is homogeneous itself provided $G/\Gamma$ does not
admit a real rank 1 $T$-invariant quotient. An immediate corollary
from our Theorem \ref{thm1} shows that Margulis conjecture fails for every
Hilbert modular space of rank $2$ (Corollary \ref{cor1}), for
instance, it fails when $G = \SL(2,\R) \times \SL(2,\R)$ and $\Gamma$ is the
diagonal imbedding of $\SL(2,\sqrt{2})$ in $G$. We apply this result to produce counter-exemples to \cite[Conjecture
1]{Margulis problems survey} for much larger classes of homogeneous
spaces as $\mathrm{SO}(f,\R)/\mathrm{SO}(f,\Z)$ and
$\SL(n,\R)/\SL(n,\Z), n \geq 4$ (see Corollary \ref{Weil} and \S7).
For actions of split tori on $\SL(n,\R)/\SL(n,\Z)$ completely
different examples of orbits with non-homogeneous closures
contradicting \cite[Conjecture 1]{Margulis problems survey} have
been first constructed by F.Maucourant \cite{Mau} when $n \geq 6$
and by U.Shapira \cite{Shapira} when $n = 3$. The constructions from
\cite{Mau} and \cite{Shapira} do not apply to the class of Hilbert
modular spaces.

It is instructive to note that the split tori action on homogeneous
spaces with finite volume is the counterpart of the unipotent
subgroups action on such spaces. The latter action is completely
understood in most general setting by M.Ratner \cite{Ratner}. (See
also the earlier papers \cite{Dani-Margulis}, \cite{Margulis
Oppenheim 2}, \cite{Shah}.) One of the basic intrinsic differences
between the two actions is that the unipotent orbits never diverge.
This is a fundamental result of Margulis \cite{Margulis
non-divergence} which admits important quantitatif versions
(cf.\cite{Dani}, \cite{Dani-Margulis+}, \cite{Kleinbock-Margulis}).

We describe in this paper the closures of locally divergent
$D_I$-orbits on the Hilbert modular spaces $G/\Gamma$. It turns out
that, on one hand, Conjecture A is not valid for the action of
two-dimensional tori (Theorem \ref{thm1}) and in the case of Hilbert
modular spaces corresponding to $\mathrm{CM}$-fields (Theorem
\ref{thm3}) and, on the other hand, Conjecture A is valid in all
remaining cases (Theorem \ref{thm2}).

Let us formulate our theorems. The cases $\#I = 2$ and $\#I
> 2$ are very different by nature and will be considered
separately.

\begin{thm}
\label{thm1} Let $\#I = 2$ and $D_I\pi(g)$ be a locally divergent
orbit on $G/\Gamma$. Suppose that the closure $\overline{D_I\pi(g)}$
is not an orbit of a torus. Then
$$\overline{D_I\pi(g)} = \underset{i=1}{\overset{s}{\bigcup}}T_i\pi(h_i) \bigcup {D_I\pi(g)},$$
where $2 \leq s \leq 4$, $T_i$ are tori containing $D_I$ and
$T_i\pi(h_i)$ are pairwise different closed non-compact orbits. In
particular, if $\#I = 2$ then there are no dense locally divergent
$D_I$-orbits.
\end{thm}

The locally divergent orbits $D_I\pi(g), \ \#I \geq 2,$ such that
$\overline{D_I\pi(g)}$ is not an orbit of a torus always exist and
are explicitly described by Corollary \ref{cor3} below. Moreover, as
shown by Proposition \ref{s=4}, there are locally divergent orbits
for which the boundaries of their closures consist of exactly $s =
4$ different closed orbits.

Theorem \ref{thm1} easily implies that the orbit rigidity conjecture
in the case of Hilbert modular spaces is not valid. More precisely,
we have the following:

\begin{cor}
\label{cor0} Let $\#I = 2$ and $T = D_I \ \mathrm{or} \ D_{I,\R}$.
Suppose that $T\pi(g)$ is a locally divergent orbit such that
$\overline{T\pi(g)}$ is not an orbit of a torus. Then the orbit
$T\pi(g)$ is a proper open subset of $\overline{T\pi(g)}$. In
particular, $\overline{T\pi(g)}$ is not homogeneous.
\end{cor}

The maximal tori action (the so-called Weyl chamber flow) deserves
special attention. The next corollary is a particular case of
Theorem \ref{thm1}:

\begin{cor}
\label{cor1} Suppose that the Hilbert modular space $G/\Gamma$ is of
rank $r = 2$. Then a locally divergent orbit $D\pi(g)$ is either
closed or
$$\overline{D\pi(g)} \setminus {D \pi(g)} =
\underset{i=1}{\overset{s}{\bigcup}}D\pi(h_i),$$ where $2 \leq s \leq
4$, and $D\pi(h_i)$ are pairwise different, closed, non-compact
orbits.
\end{cor}

After the main results of this paper had been reported \cite{Toma2},
appeared the preprint of E.lindenstrauss and U.Shapira \cite{LS}
where, using different ideas, the authors prove a somewhat similar
to the above corollary result for the action of maximal tori on
$\SL(3,\R)/\SL(3,\Z)$.

The homogeneous space $G/\Gamma$ in the formulation of Corollary
\ref{cor1} can be embedded as a closed subvariety in a number of
homogeneous spaces $H/\Delta$ where $H$ is a semi-simple Lie group
and $\Delta$ its irreducible lattice. We use this to obtain more
examples of multidimensional tori orbits with non-homogeneous
closures. For instance, at the end of Section 4 we will prove:

\begin{cor}
\label{Weil} Suppose that one of the following holds:
\begin{enumerate}
\item[(a)] $H = \mathrm{SO}(f,\R)$ and $\Delta = \mathrm{SO}(f,\Z)$, where $f$ is a non-degenerate
quadratic form with rational coefficients of $n \geq 5$ variables, of
$\R$-rank $\geq 2$, and of $\Q$-rank $\geq 1$;
\item[(b)] $H = \SL(n,\R)$, $\Delta = \SL(n,\Z)$ and $n \geq 4$.
\end{enumerate}
Let $T$ be a maximal $\R$-split torus of $H$ acting on $H/\Delta$ by
left multiplication and let $\pi_{\circ}:H \rightarrow H/\Delta, \ g
\mapsto g\Delta$. Then there exist orbits $T\pi_{\circ}(g)$ such
that
$$\overline{T\pi_{\circ}(g)} \setminus {T \pi_{\circ}(g)} =
\underset{i=1}{\overset{4}{\cup}}T\pi_{\circ}(h_i)$$ where
$T\pi_{\circ}(h_i), \ 1\leq i \leq 4,$ are pairwise different,
closed, non-compact orbits.
\end{cor}

Recall that if $f$ is a real isotropic quadratic form of $n = 3$
variables then $\mathrm{SO}(f,\R)$ is locally isomorphic to
$\SL(2,\R)$. If $n = 4$, $\rank_{\Q}f = 1$ and $\rank_{\R}f = 2$
then $\mathrm{SO}(f,\R)$ is locally isomorphic to $\SL(2,\R) \times
\SL(2,\R)$ and $\mathrm{SO}(f,\Z)$ is an irreducible non-uniform
lattice in $\mathrm{SO}(f,\R)$ (cf.\cite[Theorems 5.21 and
5.22]{Artin}). If, $n \geq 5$, then $\mathrm{SO}(f)$ is a simple group
of type $B_{\frac{n-1}{2}}$ if $n$ is odd and of type
$D_{\frac{n}{2}}$ if $n$ is even.

The dynamics of the action of $D_I$ on a Hilbert modular space
$G/\Gamma$ differs drastically when $\#I > 2$. In this case the
so-called {\it $\mathrm{CM}$-fields} play an important role. Recall
that a number field $K$ is called {\it $\mathrm{CM}$}-field (so
named for a close connection to the theory of complex
multiplication) if it is a quadratic extension of a totally real
number field which is totally imaginary.

\begin{thm}
\label{thm2} Let $\#I > 2$ and $D_I\pi(g)$ be a locally divergent
orbit such that $\overline{D_I\pi(g)}$ is not an orbit of a torus.
Assume that $K$ is not a $\mathrm{CM}$-field. Then ${D_I\pi(g)}$ is
a dense orbit.
\end{thm}

If $K$ {\it is} a $\mathrm{CM}$-field then the closure of
${D_I\pi(g)}$ might not be homogeneous. This is related to a simple
observation which we are going to explain now. Denote by $G_{i,\R},
\ 1 \leq i \leq r,$ the subgroup of real matrices in $G_i$ and put
$G_{\R} = \underset{i=1}{\overset{r}{\prod}}G_{i,\R}$. Clearly,
$G_{\R} \supset D_{I,\R}$. Now let $K$ be a $\mathrm{CM}$-field
which is a quadratic extension of a totally real number field $F$
and let $\OO_F$ be the ring of integers of $F$. Then $\Gamma_{\R} =
\SL(2,\OO_F)$ is a lattice in $G_{\R}$ and the orbit $G_{\R}\pi(e)$
is closed and homeomorphic to $G_{\R}/\Gamma_{\R}$. It is standard
to prove that this property characterize $K$ as a
$\mathrm{CM}$-field, that is, if $G/\Gamma$ admits a closed
$G_{\R}$-orbit then $K$ is a $\mathrm{CM}$-field. It follows from
the special case of Theorem \ref{thm2} for totally real fields
(Corollary \ref{real}) that if $K$ is a $\mathrm{CM}$-field of
degree $> 4$, $x \in G_{\R}\pi(e)$ and $D_{I,\R}x$ is a locally
divergent orbit whose closure is not an orbit of a torus, then
$\overline{D_{I,\R} x} = G_{\R}\pi(e)$. Since $D_I$ is an extension
of $D_{I,\R}$ by a compact torus this implies that $\overline{D_I x}
= D_IG_{\R}\pi(e)$. It is clear that $\overline{D_I x}$ is not
homogeneous which shows that if $K$ is a $\mathrm{CM}$-field the
analog of Theorem \ref{thm2} is not valid.

Let us turn to the study of the orbits for the action of the
$\R$-split tori $D_{I,\R}$ which is also important from the point of
view of Margulis' conjecture.

In the classical case of {\it real} Hilbert modular spaces in view
of Theorem \ref{thm2} we have:

\begin{cor}
\label{real} Let $K$ be a totally real number field of degree $r
\geq 3$. Let $\#I > 2$ and $D_I\pi(g)$ be a locally divergent orbit
such that $\overline{D_I\pi(g)}$ is not an orbit of a torus. Then
$\overline{D_I\pi(g)} = G/\Gamma$.

In particular, if $D_I = D$ then ${D\pi(g)}$ is either closed or
dense.
\end{cor}

In \S5 we prove the following generalization of Corollary
\ref{real}:

\begin{cor}
\label{cor2}
With the assumptions of Theorem \ref{thm2}, the orbit
$D_{I,\R}\pi(g)$ is dense in $G/\Gamma$.
\end{cor}

When $K$ is a $\mathrm{CM}$-field we obtain examples of tori orbits
contradicting Margulis' conjecture which are {\it essentially}
different from those provided by Theorem \ref{thm1}.

\begin{thm}
\label{thm3} Let $K$ be a $\mathrm{CM}$-field and $\#I
> 2$. Then there exists a point $x \in G/\Gamma$ with the following
properties:
\begin{enumerate}
\item[(i)] The orbit $D_{I,\R}x$ is locally divergent and $\overline{D_{I,\R}x} \neq G/\Gamma$;
\item[(ii)] There exists an $y \in \overline{D_{I,\R}x} \setminus
D_{I,\R}x$ such that $\overline{D_{I,\R}x} = \overline{D_{I,\R}y} $
and $Hy$ is not closed for any proper subgroup $H$ of $G$ containing
$D_{I,\R}$;
\item[(iii)] $\overline{D_{I,\R}x} \setminus D_{I,\R}x$ is not contain in a union of countably
many closed orbits of proper subgroups of $G$.

In particular, $\overline{D_{I,\R}x}$ is not homogeneous.
\end{enumerate}
\end{thm}

As a by-product of the proofs of the above theorems we get the
following corollary which is known for $D_I = D$  (see Theorem
\ref{q+} below).

\begin{cor}
\label{cor3} Let $\#I \geq 2$. Then $\mathcal{N}_{G}(D_I)G_K
\subsetneqq \bigcap_{i \in I}\big(\mathcal{N}_{G}(D_i)G_K\big)$ and
$D_I\pi(g)$ is a locally divergent orbit such that
$\overline{D_I\pi(g)}$ is not an orbit of a torus if and only
$$g \in \bigcap_{i \in I}\big(\mathcal{N}_{G}(D_i)G_K \big) \setminus \mathcal{N}_{G}(D_I)G_K.$$
\end{cor}

\medskip

The following orbit rigidity conjecture
is plausible:

\medskip

\textbf{Conjecture B.} Let $G$ be a real semisimple algebraic group
with no compact factors and let $\Gamma$ be an irreducible lattice
in $G$. Suppose that $\rank_{\R} G \geq 2$ and that every semisimple
subgroup $G_0$ in $G$ of the same $\R$-rank as $G$ acts minimally on
$G/\Gamma$ (i.e., every $G_0$-orbit is dense). Then if $T$ is a
maximal $\R$-split torus in $G$ and $x \in G/\Gamma$, either
\begin{enumerate}
\item $\overline{Tx} = G/\Gamma$, or
\item $\overline{Tx}\setminus {Tx}\subset \underset{i=1}{\overset{n}{\bigcup}} H_ix_i$, where $H_i$
are proper reductive subgroups of $G$ containing $T$ and the orbits
$H_ix_i$ are closed.
\end{enumerate}

\medskip

We apply our method to study the values of binary quadratic forms at
integral points. Denote $A =
\overset{r}{\underset{i=1}{\prod}}K_{i}$ and $A^* =
\overset{r}{\underset{i=1}{\prod}}K_{i}^*$. The polynomial ring
$A[X,Y]$ is naturally isomorphic to
$\overset{r}{\underset{i=1}{\prod}}K_{i}[X,Y]$. The natural
embeddings of $K$ into $K_i$ induce embeddings of $K[X,Y]$ into
$K_i[X,Y]$, $1 \leq i \leq r,$ and a diagonal embedding of $K[X,Y]$
into $A[X,Y]$. {\it In the next theorem $f = (f_i)_{i\in
\overline{1,r}} \in A[X,Y]$, where $f_i \in K_i[X,Y]$ are split,
non-degenerate, quadratic forms over $K$} (that is, $f_i = l_{i,1}
\cdot l_{i,2}$, where $l_{i,1}$ and $l_{i,2}$ are linearly
independent linear forms with coefficients from $K$). If
$(\alpha,\beta) \in \OO^2$ then $f(\alpha,\beta)$ is an element in
$A$ with its $i$-th coordinate equal to $f_i(\alpha,\beta)$. It is
clear that if $f_i$ are two by two proportional (equivalently, if
there exists a $g \in K[X,Y]$ such that $f_i = c_i \cdot g$, $c_i
\in K,$ for all $i$) then $f(\OO^2)$ is a discrete subset of $A$. It
follows from \cite[Theorem 1.8]{Toma} that the opposite is also
valid: the discreteness of $f(\OO^2)$ in $A$ implies the
proportionality of $f_i \ , 1 \leq i \leq r$. In the next theorem we
describe the closure of $f(\OO^2)$ in $A$ when $f_i, \ 1 \leq i \leq
r,$ are not proportional.

\begin{thm}
\label{application} With the above notation and assumptions, suppose
that $f_i$ are not proportional. Then the following assertions hold:
\begin{enumerate}
\item[(a)] If $r > 2$ and $K$ is not a $\mathrm{CM}$-field then $f(\OO^2)$ is dense in
$A$;
\item[(b)] Let $r = 2$. Put $K_1' = \{f_1(x,y): (x,y) \in K_1^2 \ \mathrm{and} \ f_2(x,y) = 0\}$ and
$K_2' = \{f_2(x,y): (x,y) \in K_2^2 \ \mathrm{and} \ f_1(x,y) =
0\}$. Then there exist  $2 \leq s \leq 4$ and pairwise nonproportional $K$-rational quadratic forms $\phi^{(j)} \in
K[X,Y]$, $1 \leq j \leq s$, such that
$$
\overline{f(\OO^2)} = \underset{j=1}{\overset{s}{\bigcup}}
\phi^{(j)}(\OO^2) \bigcup (K_1' \times \{0\}) \bigcup (\{0\} \times
K_2') \bigcup f(\OO^2).
$$
So, the set $\overline{f(\OO^2)} \cap A^*$ is countable and the set
$\overline{f(\OO^2)} \bigcap (A \setminus A^*)$ is continuum.
Moreover, $K_i' = \C$ if $K_i = \C$ and $K_i' = \R, \R_- \
\mathrm{or} \ \R_+$ if $K_i = \R$.
\end{enumerate}
\end{thm}

Let us describe the organisation and the main points of the paper.
In \S2 we recall some results from our previous paper \cite{Toma}
and we prove auxiliary results about the structure of the group of
units of a number field. The phenomenon which is at the base of the
difference between the 2-dimensional tori action (Theorem
\ref{thm1}) and the higher dimensional tori action (Theorem
\ref{thm2}) is the simple fact that the projection of the group of
units to any archimedean completion $K_v^*$ of $K^*$ is discrete if
$r = 2$ and is not discrete if $r > 2$. In \S3 we use dynamical type
arguments in combination with Minkowski's theorem for the convex
body, the structure of the locally divergent orbits \cite{Toma} and the Bruhat
decomposition for $\SL_2$ in order to describe in a very explicit
way the accumulation points of the tori orbits under consideration.
In \S4 we apply these results to deduce Theorem \ref{thm1} and its
corollaries. In \S5 we use the above mention phenomenon in order to
prove that in the case of action of tori of dimension $> 2$ if
the closure of an orbit is not contained in an orbit of a larger
torus then it contains curves which approximate arbitrary long
pieces of real unipotent orbits. This allows to prove Theorem
\ref{thm2} and its corollaries using well-known properties of
unipotent flows. The proof of Theorem \ref{thm3} is a result of a
careful analyses of the previous arguments in this section. Our
number-theoretic application is proved in \S6. \S7 contains a
specification of Theorem \ref{thm1} and indications for forthcoming
works related with the paper.

The main results of the paper have been announced in \cite{Toma2}.

\section{Preliminaries} \label{section:prelims}

\subsection{Notation}
As usual, $\Q$, $\R$, and $\C$ denote the rational, real and complex
numbers, respectively. Also, $\R_+$ (respectively, $\R_-$) is the
set of nonnegatives (respectively, nonpositives) real numbers. Let
$\R_{>0} = \R_+ \setminus \{0\}$. We denote by $|\ . \ |$ the
standard norms on $\R$ and $\C$.

In this paper ${K}$ is a number field and $K_1, \cdots, K_r$ are the
completions of $K$ with respect to the archimedean places of $K$. We
denote by $|\ \cdot \ |_{i}$ {\it the normalized} valuation on
$K_i$. So, if $x \in K$ and $K_i = \R$ (respectively, $K_i = \C$)
then $|x|_{i} = | \sigma_{i}(x)|$ (respectively, $|x |_{i} = |
\sigma_{i}(x)|^2$) where $\sigma_{i}$ is the corresponding embedding
of $K$ into $K_i$. Note that $|\mathrm{N}_{K/\Q}(x)| = |x|_{1}
\cdots |x|_{r}$, where $\mathrm{N}_{K/\Q}(x)$ is the algebraic norm
of $x$. The elements of $K$ are identified with their images in
$K_{i}$ via the embeddings $\sigma_{i}$. So, if $x \in K$, with some
abuse of notation, we write $x$ instead of $\sigma_{i}(x)$. The
exact meaning of $x$ will be always clear from the context.

If $R$ is a ring then $R^*$ is its group of invertible elements.

Let $A = \overset{r}{\underset{i=1}{\prod}}K_i$ and $A^* =
\overset{r}{\underset{i=1}{\prod}}K_{i}^*$. $A$ (respectively,
$A^*$) is a topological ring (respectively, topological group)
endowed with the product topology. The field $K$ (respectively, the
group $K^*$) is diagonally embedded in $A$ (respectively, $A^*$).
The ring of integers $\OO$ of $K$ is a co-compact lattice of $A$ and
the group of units $\OO^*$ is a discrete subgroup of $A^*$.

If $M$ is a subset of a topological space $X$ then $\overline{M}$ is
the topological closure of $M$ in $X$. Also, if $H$ is a closed
subgroup of a topological group $L$ we denote by $H^{\circ}$ the
connected component of the identity of $H$. By
$\mathcal{N}_L(H)$ (respectively, $\mathcal{Z}_L(H)$) we denote the
normalizer (respectively, the centralizer) of $H$ in $L$.

The notation $G_i$, $G$, $G_{\R}$, $D_I$, $D_{I,\R}$ have been
introduced in the Introduction. The group $G$ is considered as a
{\it real} Lie group.

The diagonal embedding of $\SL(2,K)$ in $G$ will be denoted by
$G_K$. $B_K^+, B_K^-$ and $D_K$ are the groups of upper triangular,
lower triangular and diagonal matrices in $G_K$, respectively. For
every $1 \leq i \leq r$ we denote by $G_{i,K}$, $B_{i,K}^+,
B_{i,K}^-$ and $D_{i,K}$ the images of $G_K$, $B_K^+, B_K^-$ and
$D_K$, respectively, under the natural projection $G \rightarrow
G_{i}$.

In the course of our considerations one and the same matrix with
coefficients from $K$ might be considered, according to the context,
as an element from $G_K$ or from $G_{i,K}$. For instance, if $g =
(g_1, \cdots, g_r) \in G$ and $g_i \in G_{i,K}$ writing $\pi(g_i)$,
where $\pi$ is the map $G \rightarrow G/\Gamma$, $g \mapsto
g\Gamma$, we mean that $g_i$ is considered as an element from $G$
and, therefore, from $G_K$.

Given a non-empty subset $I$ of $\{1, \cdots, r \}$ we put $A_I^*
\bydefn {\underset{i\in I}{\prod}}K_i^*$. Let
$d_i: K_i^* \rightarrow G_i, x \mapsto \left(\begin{array}{cc} x&0\\
0&x^{-1}\\ \end{array} \right)$. We put $d_I \bydefn {\underset{i\in
I}{\prod}}d_i$ and $d \bydefn d_{\{1, \cdots, r \}}$. So, $D_I =
d_I((A_I^*) ^\circ)$.

Let $\mathfrak{g}_i = \mathfrak{sl}(2,K_i)$, $\mathfrak{g} =
\overset{r}{\underset{i=1}{\prod}}\mathfrak{g}_i$, $\mathfrak{g}_K =
\mathfrak{sl}(2,K)$ and $\mathfrak{g}_{\OO} = \mathfrak{sl}(2,\OO)$.
Fixing a basis of $K$-rational vectors in $\mathfrak{g}_K$ we denote
by $\parallel \cdot \parallel_i$ the norm $\max$ on
$\mathfrak{g}_i$. Since $\mathfrak{g} =
\overset{r}{\underset{i=1}{\prod}}\mathfrak{g}_i$ we can define a
norm $\parallel \cdot \parallel$
 on $\mathfrak{g}$ by $\parallel
\mathbf{x} \parallel =
\max\limits_i\parallel\mathbf{x}_{i}\parallel_i$, $\mathbf{x}=
(\mathbf{x}_1, \cdots, \mathbf{x}_r) \in \mathfrak{g}$.

As usual, we denote by $\Ad: G \rightarrow
\mathrm{Aut}(\mathfrak{g})$ the adjoint representation of $G$.

\subsection{Locally divergent orbits}
\label{locally divergent orbits}

The following theorem is a very particular case of \cite[Theorem
1.4]{Toma} (see also \cite[Corollary 1.7]{Toma}). The paper \cite{Toma}
is related with \cite{Tomanov-Weiss}. Prior to \cite{Tomanov-Weiss}
Margulis described the divergent orbits for the action of the full
diagonal group on the space of lattices of $\R^n, \ n \geq 2$
\cite[Appendix]{Tomanov-Weiss}.

\begin{thm}
\label{q+} Let $r \geq 2$, $g = (g_1, \cdots, g_r) \in G$, and $I$
be a non-empty subset of $\{1, \cdots, r\}$. The following
assertions hold:
\begin{enumerate}
\item[(a)] If the orbit $D_I\pi(g)$ is closed then either $I$ is a singleton
or $I = \{1, \cdots, r\}$;
\item[(b)]$D_i\pi(g), 1 \leq i \leq r,$ is closed (equivalently, divergent) if and only if $g \in \mathcal{N}_{G}(D_i)G_K$
(equivalently, $g_i \in D_iG_{i,K}$);
\item[(c)] The following conditions are equivalent:
\begin{enumerate}
\item [(i)]$D\pi(g)$ is closed and non-compact;
\item [(ii)] $D\pi(g)$ is closed and locally divergent;
\item [(iii)] $g \in \mathcal{N}_{G}(D)G_K$.
\end{enumerate}
\end{enumerate}
\end{thm}

We will need the following proposition:

\begin{prop}
\label{lem} If $g \in \mathcal{N}_G(D_I)G_K$ then
$\overline{D_I\pi(g)} = T\pi(g)$ where $T$ is a torus containing
$D_I$.
\end{prop}
{\bf Proof.} In view of our assumption $g = g'h$ where $h \in
\mathcal{N}_G(D)G_K$ and $g' \in \underset{i \notin I}\prod G_{i}$.
Let $\Delta$ be the stabilizer of $\pi(g)$ in $g'Dg'^{-1}$. It
follows from Theorem \ref{q+}(c) that $g'D\pi(h)$ is closed. Since
$\overline{D_I\pi(g)} \subset g'Dg'^{-1}\pi(g)$ we get that
$\overline{D_I\pi(g)} = T\pi(g)$ where $T$ is the connected
component of the identity of $\overline{D_I \Delta}$. \qed

\subsection{Propositions about the units}

Denote $A^1 = \{(x_1, \cdots, x_r) \in A^*: |x_1|_1 \cdots |x_r|_r =
1\}$. Given a positive integer $m$ we put $\OO^*_m = \{\xi^m | \xi
\in \OO^*\}$.

The following lemma follows easily from the classical fact that
$\OO^*$ is a lattice in $A^1$.

\begin{lem}
\label{lem0}\mbox{\rm{(}}cf.\cite[Lemma 3.2]{Toma}\mbox{\rm{)}} Let
$m$ be a positive integer. There exists a real $\kappa_m
> 1$ with the following property. Let  ${x} =
({x}_{i}) \in A^*$ and for each $1 \leq i \leq r$ let $a_i$ be a
positive real number such that ${\underset{i}{\prod}}a_i =
{\underset{i}{\prod}}|x_i|_i$. Then there exists $\xi \in
 \OO^*_m$ such that
$$
\frac{a_i}{\kappa_m} \leq |\xi {x}_{i}|_i \leq \kappa_m a_i
$$
for all $i$.
\end{lem}

\begin{prop}
\label{integers1} Let $r \geq 3$, $3 \leq l \leq r$, $I = \{l,
\cdots, r\}$ and $p_I: A^* \rightarrow A_I^*$ be the natural
projection. Denote by $H$ the closure of $p_I(\OO^*)$ in $A_{I}^*$.
Then
\begin{enumerate}
\item[(a)] the projection of $H^\circ$ into each $K_i^*, i \geq l,$ is non-trivial;
\item[(b)] for any real $C > 1$ there exists $\xi \in \OO^*$ such that $|\xi|_l > C$
and $|1-|\xi|_i| < \frac{1}{C}$ for all $i > l$.
\end{enumerate}
\end{prop}
\medskip
{\bf Proof.} (a) By Dirichlet's theorem for the units there exists a
positive integer $m$ such that $\OO^*_m$ is a free abelian group of
rank $r-1$. It is clear that $H^\circ$ coincides with the connected
component of the closure of $p_I(\OO_m^*)$. Since $H^{\circ}$ is
open in $H$ and $\OO^*_m$ is diagonally embedded in $H$ it is enough
to show that $H^\circ \neq \{1\}$. Suppose that $H^\circ =
\{1\}$. Then $H$ is a discrete subgroup of $A_I^*$ containing a free
subgroup of rank $r-1$. This is a contradiction because $l \geq 3$ and $A_I^*$ is a
direct product of a compact group and $\Z^{r-l +1}$.

\medskip

(b) Consider the logarithmic representation of the group of units
$\log_S: \OO^* \rightarrow \R^r, \theta \mapsto (\log|\theta|_1,
\cdots, \log|\theta|_r)$ (see \cite{W}). According to the Dirichlet
theorem $\log_S(\OO^*)$ is a lattice in the hyperplane $L = \{(x_1,
\cdots, x_r) \in \R^r : x_1 + x_2 + \cdots + x_r = 0\}$. Let $\psi:
L \rightarrow \R^{r-1}, (x_1, \cdots, x_r)$ $\mapsto (x_2, \cdots,
x_r)$. Then $\psi(\log_S(\OO^*))$ is a lattice in $\R^{r-1}$ with
co-volume equal to a positive real $V$. For every natural $n$, we denote
$$B_n = \{(x_2, \cdots, x_r) \in \R^{r-1}: |x_i| \leq \frac{1}{n} \
\mathrm{if} \ i \neq l
 \ \mathrm{and} \ |x_l| \leq {n}^{r-2}V\}.$$
By Minkowski's lemma there exists a $\xi_n \in \OO^*$ such that
$\psi(\log_S(\xi_n)) \in B_n \setminus \{0\}$. If the sequence
$|\xi_n|_l$ is unbounded from above then we can choose $\xi = \xi_n$
with $n$ large enough. Let $|\xi_n|_l < C$ where $C$ is a constant.
Since $\psi(\log_S(\OO^*))$ is discrete this implies the existence
of a unit $\eta$ of infinite order such that $|\eta|_l
> 1$ and $|\eta|_i = 1$ if $i \neq l$ and $i > 1$. Hence we can choose $\xi =
\eta^m$ with $m$ sufficiently large. \qed

\medskip

\begin{prop}
\label{integers2} Let $p_l: A^* \rightarrow K_{l}^*$, $1 \leq l \leq
r$, be the natural projection. Assume that $K_l = \C$ and that the
connected component of the identity of $\overline{p_l(\OO^*)}$ coincides with
$\R_{>0}$. Then $K$ is a $\mathrm{CM}$-field.
\end{prop}

\medskip

{\bf Proof.} There exists a positive integer $m$ such that
$\overline{p_l(\OO_m^*)} = \R_{>0}$. Denote by $F$ the subfield of
$K$ generated over $\Q$ by all $\theta \in \OO_m^*$ and denote by
$\OO_F^*$ the group of units of $F$. Let $s$, respectively $t$, be
the number of real, respectively complex, places of $K$ and let
$s_1$, respectively $t_1$, be the number of real, respectively
complex, places of $F$. By Dirichlet's theorem $\OO_m^*$ is a free
group of rank $s+t-1$. Since $\OO_m^* \subset \OO_F^* \subset \OO^*$
and the group of principal units of $F$ is free of rank $s_1+t_1-1$
we have
$$r-1 = s+t-1 = s_1+t_1-1.$$
Let $n$ be the degree of $K$ over $F$. Since $s+2t$ is the
degree of $K$ over $\Q$ and $s_1+2t_1$ is the degree of $F$ over
$\Q$ we get
\begin{eqnarray}
s + 2t = n(s_1 + 2t_1) \Leftrightarrow r+t = n(r+t_1)
\Leftrightarrow \nonumber \\
(n-1)r  = t-t_1n \Leftrightarrow (n-1)(t+s) = t-t_1n. \nonumber
\end{eqnarray}
Since $n>1$ the last equality implies that $s = t_1 =0$ and $n=2$
proving the proposition. \qed

\medskip

\textbf{Example.}{\footnote{This
example is essentially due to Yves Benoist.}} There are non-$\mathrm{CM}$ fields such that the
connected components of the identity of $\overline{{p_l(\OO^*)}}$ are 1-dimensional subgroups
of $\C^*$ different from $\R_{>0}$. Such fields need
special treatment in the course of the proof of Proposition
\ref{prop4}(a) below. An example of this type is provided by the
field $K = \Q(\alpha)$ where $\alpha$ is a root of the equation $(x
+ \frac{1}{x})^2 - 2(x + \frac{1}{x}) - 1 = 0$. The field $K$ has two
real and one (up to conjugation) complex completions. If $K_3 = \C$
then it is easy to see that $\overline{p_3(\OO^*)}^\circ$ coincides
with the unit circle.

\medskip

\section{Accumulations points for locally divergent orbits}
Up to the end of the paper $D_I\pi(g)$ will denote a locally
divergent orbit. In view of Theorem \ref{q+}(b), {\it we may (and
will) assume without loss of generality that $g = (g_1, \cdots,
g_r)$ with $g_i \in G_{i,K}$ whenever $i \in I$}.

\medskip
The following lemma is an easy consequence from the commensurability
of $\Gamma$ and $h\Gamma h^{-1}$ when $h \in G_K$.

\begin{lem}
\label{lem2} Let $h \in G_K$. The following assertions hold:
\begin{enumerate}
\item[(a)] There exists a positive integer $m$
such that $d(\xi)\pi(h) = \pi(h)$ for all $\xi \in \OO^*_m$;
\item[(b)] If $\{\pi(g_i)\}$ is a
converging sequence in $G/\Gamma$ then there exists a converging
subsequence of $\{\pi(g_ih)\}$;
\item[(c)] If $\overline{D_I\pi(g)} = G/\Gamma$ then $\overline{D_I\pi(gh)} =
G/\Gamma$.
\end{enumerate}
\end{lem}

\begin{prop}\label{prop1}
 Let $I = \{1,2\}$ and $(s_k, t_k) \in K_1^* \times K_2^*$ be a sequence such that
 $|\log |s_k|_1| + |\log
|t_k|_2| \underset{k}{\rightarrow} \infty$ and $d_I(s_k,t_k)\pi(g)$
converges to an element from $G/\Gamma$. Then:
\begin{enumerate}
\item[(a)] There exists a constant $C > 1$ such that $-C < |\log |s_k|_1| -
|\log |t_k|_2 | < C$;
\item[(b)] Let $|s_k|_1 \rightarrow \infty$, $|t_k|_2 \rightarrow 0$. 
 Then
$g_1g_2^{-1} = b_-b_+^{-1}$, where $b_- \in B^-_K$ and $b_+ \in
B^+_K$.
\end{enumerate}
\end{prop}
{\bf Proof.}(a) The remaining cases being analogous, it is enough to
consider the case when $|s_k|_1 \rightarrow  \infty$ and
$\underset{k}{\sup}\frac{\max\{ |t_k|_2, |t_k|_2^{-1}\}}{|s_k|_1} < \infty$.

Assume on the contrary that (a) is false. Then $\frac{\max\{
|t_k|_2, |t_k|_2^{-1}\}}{|s_k|_1} \underset{k}{\rightarrow} 0$. It
is well known that for every $h \in G_K$ $\Ad(h)\mathfrak{g}_{\OO}$
is commensurable with $\mathfrak{g}_{\OO}$. Since $g_1 \in G_K$ this
implies the existence of $\mathbf{u} \in \Ad(g)\mathfrak{g}_{\OO},
\mathbf{u} \neq 0,$ such that $\mathrm{pr}_1(\mathbf{u})$ is a lower
triangular nilpotent matrix where $\mathrm{pr}_1$ is the projection
of $\mathfrak{g}$ to $\mathfrak{g}_1$. Recall that $\mathfrak{g} =
\overset{r}{\underset{i=1}{\prod}}\mathfrak{g}_i$. Let
$\Ad(d_I(s_k, t_k))(\mathbf{u}) = (\mathbf{u}_1^{(k)}, \cdots,
\mathbf{u}_r^{(k)}) \in \mathfrak{g}$. Since $\frac{\max\{ |t_k|_2,
|t_k|_2^{-1}\}}{|s_k|_1} \underset{k}{\rightarrow} 0$ and
$\Ad(d_I(s_k, t_k))$ is acting by conjugation on the elements from
$\mathfrak{g}$, we see that $\|\mathbf{u}_1^{(k)}\|_1 \cdots
\|\mathbf{u}_r^{(k)}\|_r \underset{k}{\rightarrow} 0$.
In view of Lemma \ref{lem0}, there exists a sequence $\xi_k \in
\OO^*$ such that $\parallel \Ad(d_I(s_k, t_k))(\xi_k\mathbf{u})
\parallel = \|(\xi_k\mathbf{u}_1^{(k)}, \cdots, \xi_k\mathbf{u}_r^{(k)})\| \underset{k}{\rightarrow} 0$.
By Mahler's compactness criterion $d_I(s_k,t_k)\pi(g)$ tends to
infinity which is a contradiction.

\medskip

(b) By Bruhat decomposition
$$
G_K = B^+_K\cup B^+_K\omega B^+_K = \omega B^+_K \cup B^-_KB^+_K,
$$
where $\omega = \left(
\begin{array}{cc} 0&1\\ -1&0\\ \end{array} \right)$.

Suppose on the contrary that $g_1g_2^{-1} \in \omega B^+_K$.
Shifting $g$ from the right by $g_2^{-1}$ and from the left by a
suitable element from $\mathcal{Z}_G(D_I)$ we may (and will) assume
with no loss of generality (see Lemma \ref{lem2}(b)) that $g_1 =
\omega u^+(\alpha)$, where $u^+(\alpha) = \left(
\begin{array}{cc} 1&\alpha\\ 0&1\\ \end{array} \right)$, $\alpha \in K$, and $g_i = e$
for all $i > 1$. In view of (a), there exists a constant $C > 1$
such that $\frac{1}{C} < |s_k|_1 \cdot |t_k|_2 < C$. Now using Lemma
\ref{lem2}(a) and Lemma \ref{lem0} we find a sequence $\xi_k \in
\OO^*$ and a positif constant $\kappa$ such that
$d(\xi_k)\pi(u^+(\alpha)) = \pi(u^+(\alpha))$, $\frac{1}{\kappa} <
\frac{|s_k|_1}{|\xi_k|_1} < \kappa$, $\frac{1}{\kappa} <
\frac{|t_k|_2}{|\xi_k|_2} < \kappa$ and $\frac{1}{\kappa} <
|\xi_k|_i < \kappa$ for all $i > 2$. Then $(s_k,t_k, 1, \cdots,1) =
\xi_k a_k$ where $a_k \in A^*$ is a bounded sequence. Passing to a
subsequence we can suppose that $a_k$ converges to an element from
$A^*$. Note that $d(\xi_k)\pi(g)$ converges in $G/\Gamma$ because
$d_I(s_k,t_k)\pi(g)$ does.

By an easy computation:
\begin{eqnarray*}
& & d(\xi_k)\pi(g) = d(\xi_k)(\omega u^+(\alpha), e, \cdots,
e)\pi(e) = \\ & & d(\xi_k)(\omega, u^+(-\alpha), \cdots,
u^+(-\alpha))\pi(u^+(\alpha)) = \\ & & d(\xi_k)(\omega,
u^+(-\alpha), \cdots, u^+(-\alpha))d(\xi_k^{-1})\pi(u^+(\alpha)) = \\
& & (\omega, u^+(-\alpha \xi_k^2), \cdots, u^+(-\alpha
\xi_k^2))(d_1(\xi_k^{-2}), e, \cdots, e)\pi(u^+(\alpha)).
\end{eqnarray*}

In view of the choice of $\xi_k$ we have that $|\xi_k|_1 \rightarrow
\infty$ and $|\xi_k|_2 \rightarrow 0$. Hence $|\xi_k|_i < \kappa$ if
$i \geq 2$ and $k$ is sufficiently large. So, after passing to a
subsequence, $(\omega, u^+(-\alpha \xi_k^2), \cdots, \\ u^+(-\alpha
\xi_k^2))$ converges to an element from $G$ and $d_1(\xi_k^{-2})$
tends to infinity. The latter contradicts the convergence of the
sequence $d(\xi_k)\pi(g)$. Therefore, $g_1g_2^{-1} \in B^-_KB^+_K$.
\qed

\medskip

\begin{prop}\label{prop2}
 Let $I = \{1, \cdots, l\}$ where $1 < l \leq r$, $g_1 = \cdots = g_{l-1}$ and $g_1g_l^{-1} = b_-b_+^{-1}$
where $b_- \in B^-_K$ and $b_+ \in B^+_K$. Denote $h = b_-^{-1}g_1 =
b_+^{-1}g_l$. Then we have the following:
\begin{enumerate}
\item[(a)] $(h, \cdots, h, g_{l+1}, \cdots, g_r)\pi(e) \in \overline{D_I
\pi(g)}$;
\item[(b)] Let $s_k = (s_k^{(1)}, \cdots, s_k^{(l)}) \in A_I^*$ be such
that $|s_k^{(i)}|_i \underset{k}{\rightarrow} \infty$ for all $1
\leq i < l$, $|s_k^{(l)}|_l \underset{k}{\rightarrow} 0$ and
$\frac{1}{C} < |s_k^{(1)}|_1 \cdots |s_k^{(l)}|_l < C$, where $C$ is
a positive constant. Then $d_I(s_k)\pi(g)$ admits a converging
subsequence and the limit of every such subsequence belongs to
$\overline{D_I\pi((h, \cdots, h, g_{l+1}, \cdots, g_r))}$.
\end{enumerate}
\end{prop}
{\bf Proof.} Fix $m$ such that $d(\xi)\pi(h) = \pi(h)$ for all $\xi
\in \OO^*_m$. With $s_k$ as in the formulation of (b), in view of
Lemma \ref{lem0} there exists a sequence $\xi_k \in \OO^*_m$ and a
constant $C_1 > 1$ such that $\frac{1}{C_1} <
|s_k^{(i)}\xi_k^{-1}|_i < C_1$ if $1 \leq i \leq l$ and
$\frac{1}{C_1} < |\xi_k|_i < C_1$ if $i > l$. Put $a_k =
(\underbrace{\xi_k, \cdots, \xi_k}_{l}, \underbrace{1, \cdots,
1}_{r-l})$ and $a_k' = (\underbrace{1, \cdots,
1}_{l},\underbrace{\xi_k, \cdots, \xi_k}_{r-l})$. Passing to a
subsequence we may assume that $a_k' \rightarrow a'$ where $a' \in
A^*$. In view of the choice of $\xi_k$ and the proposition
hypothesis, we get
$$\underset{k}{\lim} \ d_i(\xi_k)b_-d_i(\xi_k)^{-1} = t_- ,\ \forall \ 1 \leq i <
l,$$ and
$$
\underset{k}{\lim} \ d_l(\xi_k)b_+d_l(\xi_k)^{-1} = t_+,
$$
where $t_-$ and $t_+ \in D_K$. It is enough to prove (b) in the
particular case when $s_k^{(i)} = t_-^{-1}\xi_k$, $1 \leq i < l$,
and $s_k^{(l)} = t_+^{-1}\xi_k$.


Since $d(\xi_k)\pi(h) = \pi(h)$, we get
\begin{eqnarray*}
d_I(s_k)\pi(g) = (d_1(t_-^{-1}\xi_k)b_-, \cdots,
d_l(t_+^{-1}\xi_k)b_+, g_{l+1}h^{-1}, \cdots, g_rh^{-1})\pi(h) = \\
(d_1(t_-^{-1}\xi_k)b_-, \cdots, d_l(t_+^{-1}\xi_k)b_+,
g_{l+1}h^{-1}, \cdots, g_rh^{-1}))d({a_k}^{-1})d({a'_k}^{-1})\pi(h).
\end{eqnarray*}
Therefore
\begin{equation}
\label{one} \underset{k}{\lim} \ d_I(s_k)\pi(g) = (e, \cdots,
g_{l+1}h^{-1}, \cdots, g_rh^{-1})d(a'^{-1})\pi(h) \in \overline{D_I
\pi(g )}.
\end{equation}
Since
$$
d(a_k)^{-1}\pi(h) = d(a_k)^{-1}d(\xi_k)\pi(h) = d(a'_k)\pi(h)
\rightarrow d(a')\pi(h),
$$
multiplying (\ref{one}) by $d(a_k)^{-1}$ and passing to a limit, we
obtain that
$$
(h, \cdots, h, g_{l+1}, \cdots, g_r)\pi(e) \in \overline{D_I
\pi(g)}.
$$
Since a sequence $s_k$ with properties as in the formulation of (b)
always exists, the above proves (a). In order to complete the proof
of (b) it remains to note that
$$
\underset{k}{\lim} \ d_I(s_k)\pi(g) = \underset{k}{\lim} \
d(a_k)\pi((h, \cdots, h, g_{l+1}, \cdots, g_r)).
$$ \qed

\medskip

Let $h \in G_K$. A pair $(\sigma_1,\sigma_2) \in \{0,1\}^2$ is
called {\it admissible with respect to $h$} if
$\omega^{\sigma_1}h\omega^{\sigma_2} \in B_K^-B_K^+$, where $\omega
= \left(
\begin{array}{cc} 0&1\\ -1&0\\ \end{array} \right)$. The following lemma can be
proved by a simple calculation.

\begin{lem}
\label{admis} With the above notation, $(\sigma_1,\sigma_2)$ is
admissible with respect to  $h = \left(\begin{array}{cc} m_{11}&m_{12}\\ m_{21}&m_{22}\\
\end{array} \right)$ if
and only if $m_{1+\sigma_1,1+\sigma_2} \neq 0$.
\end{lem}

It is clear that $h \in \mathcal{N}_{G_K}(D_K)$ if and only if the
number of admissible pairs is equal to 2.

\begin{prop} \label{prop3} Let $I = \{1, \cdots, l\}$, where $1 < l < r$, $g_1 = \cdots = g_{l-1}$
and $g_1g_l^{-1} \notin \mathcal{N}_{G_K}(D_K)$. Then
$\overline{D_I\pi(g)}$ contains a point
$$(\underbrace{nh, \cdots,
nh}_{l-1},h,g_{l+1}, \cdots, g_r)\pi(e),$$ where $n \in
\mathcal{N}_{G_K}(D_K), h \in G_K$ and $hg_{l+1}^{-1} \notin
\mathcal{N}_{G_K}(D_K)$.
\end{prop}
{\bf Proof.} If the pair $(\sigma_1,\sigma_2)$ is admissible with
respect to $g_1g_l^{-1}$ then
$\omega^{\sigma_1}g_1(\omega^{\sigma_2}g_{l})^{-1} = b_-b_+^{-1}$,
where $b_- \in B^-_K$ and $b_+ \in B^+_K$, and we put
$h_{\sigma_1,\sigma_2} = b_-^{-1}\omega^{\sigma_1}g_1 =
b_+^{-1}\omega^{\sigma_2}g_l$. Shifting $\pi(g)$ from the left by
$$(\underbrace{\omega^{\sigma_1}, \cdots,
\omega^{\sigma_1}}_{l-1},\omega^{\sigma_2}, e, \cdots, e)$$ and
applying Proposition \ref{prop2}(a) we get that
$$(\underbrace{\omega^{\sigma_1}h_{\sigma_1,\sigma_2}, \cdots,
\omega^{\sigma_1}h_{\sigma_1,\sigma_2}}_{l-1},\omega^{\sigma_2}h_{\sigma_1,\sigma_2},g_{l+1},
\cdots, g_r)\pi(e) \in \overline{D_I\pi(g)}.$$

It remains to prove that $(\sigma_1,\sigma_2)$ can be chosen in such
a way that $h_{\sigma_1,\sigma_2}g_{l+1}^{-1} \notin
\mathcal{N}_{G_K}(D_K)$. Since $g_1g_l^{-1} \notin
\mathcal{N}_{G_K}(D_K)$, in view of Lemma \ref{admis} there are at
least 3 admissible pairs with respect to $g_1g_l^{-1}$. Shifting $g$
from the left by an appropriate element from
$\mathcal{N}_{G_K}(D_K)$, we may assume that $(0,0)$ and $(0,1)$ are
admissible pairs. Then
$$h_{0,0} = b_-'^{-1}g_1 = b_+'^{-1}g_2 \ \text{and} \
h_{1,0} = \widetilde{b}_-^{-1}\omega g_1 =
\widetilde{b}_+^{-1}g_2,$$ where $b_-', \widetilde{b}_- \in B^-_K$
and $b_+', \widetilde{b}_+ \in B^+_K$. Suppose on the contrary that
both $h_{0,0}g_{l+1}^{-1}$ and $h_{1,0}g_{l+1}^{-1} \in
\mathcal{N}_{G_K}(D_K)$. In view of the above expressions for
$h_{0,0}$ and $h_{1,0}$, we obtain
$$
h_{0,0}h_{1,0}^{-1} \in \mathcal{N}_{G_K}(D_K) \cap B^+_K \cap B^-_K
\omega B^-_K.
$$
This is a contradiction because $\mathcal{N}_{G_K}(D_K) \cap B^+_K =
D_K$ and $D_K \cap B^-_K \omega B^-_K = \emptyset$. \qed

\vspace{.3cm}

\section{Proofs of Theorem \ref{thm1} and Corollaries \ref{cor0} and \ref{Weil}}
\label{proofs}

\subsection{Proof of Theorem \ref{thm1}.} We suppose that $I = \{1,2\}$.
It follows from Proposition \ref{lem} that $g_1g_2^{-1} \notin
\mathcal{N}_{G_K}(D_K)$. Let $(s_k, t_k) \in K_1^* \times K_2^*$ be
an unbounded sequence such that $d_I(s_k,t_k)\pi(g)$ converges. In
view of Proposition \ref{prop1}(a) there exists a positive constant
$C$ such that $-C < |\log |s_k|_1| - |\log |t_k|_2 | < C$. Passing
to a subsequence there exist $\sigma_1$ and $\sigma_2 \in \{0, 1\}$
such that $\omega^{\sigma_1}d_1(s_k)\omega^{-\sigma_1} = d_1(s_k')$,
$\omega^{\sigma_2}d_2(t_k)\omega^{-\sigma_2} = d_2(t_k')$ where
$|s_k'|_1 \rightarrow \infty$ and $|t_k'|_2 \rightarrow 0$. Let $g'
= (\omega^{\sigma_1}g_1, \omega^{\sigma_2}g_2, g_3, \cdots, g_r)$.

It follows from Proposition \ref{prop1}(b) that
$\omega^{\sigma_1}g_1(\omega^{\sigma_2}g_2)^{-1} = b_-b_+^{-1} \in
B_K^-B_K^+$, i.e., $(\sigma_1,\sigma_2)$ is an admissible pair with
respect to $g_1g_2^{-1}$. Let
\begin{equation}
\label{boundary}
h_{\sigma_1,\sigma_2} =
b_-^{-1}\omega^{\sigma_1}g_1 = b_+^{-1}\omega^{\sigma_2}g_2.
\end{equation}
Using Proposition \ref{prop2}(b) we get:
$$
\underset{k}{\lim} \ d_I(s_k',t_k')\pi(g') \in \overline{D_I
\pi((h_{\sigma_1,\sigma_2},h_{\sigma_1,\sigma_2},g_3,\cdots,g_r))}.
$$
Therefore
$$ \underset{k}{\lim} \ d_I(s_k,t_k)\pi(g) \in \overline{D_I
\pi((\omega^{\sigma_1}h_{\sigma_1,\sigma_2},\omega^{\sigma_2}h_{\sigma_1,\sigma_2},g_3,\cdots,g_r))}.
$$
It follows that
$$
\overline{D_I\pi(g)} \subset D_I\pi(g)\cup
\cup_{(\sigma_1,\sigma_2)\in M}\overline{D_I
\pi((\omega^{\sigma_1}h_{\sigma_1,\sigma_2},\omega^{\sigma_2}h_{\sigma_1,\sigma_2},g_3,\cdots,g_r))},
$$
where $M$ is the set of all admissible pairs with respect to
$g_1g_2^{-1}$. On the other hand, using Proposition \ref{prop2}(a)
we get:
\begin{equation}
\label{limit2} \overline{D_I\pi(g)} = D_I\pi(g)\cup
\cup_{(\sigma_1,\sigma_2)\in M}\overline{D_I
\pi((\omega^{\sigma_1}h_{\sigma_1,\sigma_2},\omega^{\sigma_2}h_{\sigma_1,\sigma_2},g_3,\cdots,g_r))}.
\end{equation}
Note that
\begin{eqnarray*}
& &
\overline{D_I\pi((\omega^{\sigma_1}h_{\sigma_1,\sigma_2},\omega^{\sigma_2}h_{\sigma_1,\sigma_2},g_3,\cdots,g_r))}
= \\
& &(\omega^{\sigma_1},\omega^{\sigma_2},
g_3h_{\sigma_1,\sigma_2}^{-1},\cdots,g_rh_{\sigma_1,\sigma_2}^{-1})\overline{D_I\pi(h_{\sigma_1,\sigma_2})}.
\end{eqnarray*}
Since $D\pi(h_{\sigma_1,\sigma_2})$ is a closed locally divergent
orbit, each of the closures in the right hand side of (\ref{limit2})
is a non-compact orbit of a torus containing $D_I$. It remain to see
that at least two of these orbits are different.

Since $g_1g_2^{-1} \notin \mathcal{N}_{G_K}(D_K)$ there exists
$\sigma \in \{0,1\}$ such that $(\sigma,0)$ and $(\sigma,1) \in M$.
Suppose on the contrary that
$$
\overline{D_I\pi(\omega^{\sigma}h_{\sigma,0},h_{\sigma,0},g_3,\cdots,
g_r)} = \overline{D_I\pi(\omega^{\sigma}h_{\sigma,1},\omega
h_{\sigma,1},g_3,\cdots, g_r)}.
$$
There exist tori $T$ and $T'$ containing $D_I$ such that
$$
T\pi((h_{\sigma,0},h_{\sigma,0},g_3,\cdots,g_r)) =
T'\pi((h_{\sigma,1},\omega h_{\sigma,1},g_3,\cdots,g_r)).
$$
Then
$$
h_{\sigma,0} = th_{\sigma,1}\gamma = t' \omega h_{\sigma,1}\gamma,
$$
where $t, t' \in D_K$ and $\gamma \in \Gamma$. Hence
$$
\omega = t't^{-1}
$$
which is a contradiction. \qed

\subsection{Proof of Corollary \ref{cor0}.} We use the
notation from the formulations of Theorem \ref{thm1} and the Corollary. Let us show that both $D_I\pi(g)$
and $D_{I,\R}\pi(g)$ are open and proper in their closures. Note that if
$D_I\pi(g) \cap T_i\pi(h_i) \neq \emptyset$ for some $1 \leq i \leq
s$ then $\overline{D_I\pi(g)} \subset T_i\pi(h_i)$ which contradicts
the fact that $s \geq 2$. Therefore, the orbit $D_I\pi(g)$ is open
and proper in its closure. Suppose that there exists $i$ such that
$\overline{D_{I,\R}\pi(g)} \cap T_i\pi(h_i) = \emptyset$. Since $T_i
\supset D_I$ this implies that $\overline{D_{I}\pi(g)} \cap
T_i\pi(h_i) = \emptyset$ which is a contradiction. Therefore,
$\overline{D_{I,\R}\pi(g)} \cap T_i\pi(h_i) \neq \emptyset$ for
every $1 \leq i \leq s$. So, the orbit $D_{I,\R}\pi(g)$ is open and
proper in its closure too. Now if, supposing the contrary,
$\overline{T\pi(g)} = H\pi(g)$ for some closed subgroup $H$ then $H$
is locally homeomorphic to $T$. Since $T$ is generated by any
neighborhood of the identity, $T\pi(g)$ must be closed. This is a
contradiction completing our proof. \qed

\subsection{Proof of Corollary \ref{Weil}.}(a) It is enough to show
that $f$ represents over $\Q$ a quadratic form $f_1$ of 4 variables
such that $\rank_{\Q}f_1 = 1$ and $\rank_{\R}f_1 = 2$. Indeed, in
this case we may suppose without loss of generality that $f = f_1 +
f_2$ where $f_2$ is a quadratic form over the rationals of $n - 4$
variables. Remark that $\mathrm{SO}(f_1,\R) \times
\mathrm{SO}(f_2,\R)$ is a $\Q$-subgroup of $\mathrm{SO}(f,\R)$ and
$(\mathrm{SO}(f_1,\R) \times \mathrm{SO}(f_2,\R)) \cap
\mathrm{SO}(f,\Z)$ is commensurable with $\mathrm{SO}(f_1,\Z) \times
\mathrm{SO}(f_2,\Z)$. It is known that $\mathrm{SO}(f_1,\R) \cong
\mathrm{PSL}(2,\R) \times \mathrm{PSL}(2,\R)$ and that
$\mathrm{SO}(f_1,\Z)$ corresponds under this isomorphisme to the
diagonal embedding of $\mathrm{PSL}(2,\Z[\sqrt{d}])$ into
$\mathrm{PSL}(2,\R) \times \mathrm{PSL}(2,\R)$, where $d$ is the
discriminant of $f_1$ \cite[Theorems 5.21 and 5.22]{Artin}. If $T_1$
is a maximal $\R$-split torus of $\mathrm{SO}(f_1,\R)$ and $T_2$ is
a maximal $\R$-split torus of $\mathrm{SO}(f_2,\R)$ then $T = T_1
\times T_2$ is a maximal $\R$-split torus of $\mathrm{SO}(f,\R)$.
Now, if we choose $g_1 \in \mathrm{SO}(f_1,\R)$ in such a way that
the boundary of the closure of the orbit $T_1 g_1
\mathrm{SO}(f_1,\Z)$ consists of 4 different $T_1$-orbits
(Proposition \ref{s=4}) and if we choose $g_2 \in
\mathrm{SO}(f_2,\R)$ in such a way that the orbit $T_2 g_2
\mathrm{SO}(f_2,\Z)$ is closed (see, for example, \cite[Proposition
4.2]{Toma}) then $\overline{T\pi_\circ(g)}$, where $g = (g_1, g_2)$,
is as required.

Let us prove that $f$ represents over $\Q$ a quadratic form $f_1$
with the above mentioned properties. Since $\rank_{\Q}f \geq 1$ and
$\rank_{\R}f \geq 2$ the form $f$ is $\Q$-equivalent to a form
$x_{1}x_{2} + x_3^2 - ax_4^2 - bx_5^2 + f'(x_6, \cdots, x_{n})$
where $a$ and $b$ are rational numbers such that $a\cdot b
\neq 0$ and $b > 0$ (see \cite{Cassels}). If $b \notin \Q^2$ then we
can choose $f_1 = x_{1}x_{2} + x_3^2 - bx_5^2$. Suppose that $b \in
\Q^2$. Then the form $x_3^2 - bx_5^2$ represents a rational number
$\alpha$ such that $a \cdot \alpha \notin \Q^2$ and $a \cdot \alpha
> 0$. Therefore $f$ represents a form $f_1$ which is $\Q$-equivalent
to $x_{1}x_{2} + \alpha x_3^2 - a x_4^2$.

(b)Let $G$ and $\Gamma$ be as in the formulation of Corollary
\ref{cor1} with $K$ a real quadratic number field. Using Weil's
restriction of scalars \cite[Ch.6]{Z}, we get an injective
homomorphism $\mathrm{R}_{K/\Q}: G \rightarrow \SL(4,\R)$ such that
$\mathrm{R}_{K/\Q}(\Gamma) = \mathrm{R}_{K/\Q}(G) \cap \SL(4,\Z)$.
Let $\phi: G \rightarrow \SL(n, \R), g \mapsto \left(
\begin{array}{cc} \mathrm{R}_{K/\Q}(g)&0\\ 0&\mathrm{I}_{n-4}\\ \end{array}
\right)$, where $\mathrm{I}_{n-4}$ is the identity matrix of rank
$n-4$. Further on we identify $G$, $D$ and $\Gamma$ with $\phi(G)$,
$\phi(D)$ and $\phi(\Gamma)$, respectively. Let $F$ be the connected
component of the identity of the centralizer of $G$ in $\SL(n, \R)$.
It is clear that $F$ is a real reductive $\Q$-group, $G \cap F$ is
finite and $L = G F$ is a reductive group of real rank $n-1$. Put
$\Gamma_F = F \cap \SL(n,\Z)$. Since $L$ is a reductive $\Q$-group
the orbit $L\Gamma$ is closed in $\SL(n,\R)/\SL(n,\Z)$
(\cite[Proposition 4.2]{Toma}). Therefore the map $G/\Gamma \times
F/\Gamma_F \rightarrow \SL(n,\R)/\SL(n,\Z), \ (g\Gamma, h\Gamma_F )
\mapsto \pi_{\circ}(gh),$ is proper with finite fibers. Let $T_F$ be
a maximal $\R$-split torus in $F$ and $h \in F$ be such that $T_Fh
\Gamma_F$ is dense in $F$. Choose $g \in G$ such that the boundary
of $D\pi(g)$ consists of four pairwise different closed $D$-orbits
(Proposition \ref{s=4}). Denote $T' = DT_F$. It follows from the
above that the boundary of $T'\pi_{\circ}(gh)$ consists of four
pairwise different closed $T'$-orbits. In order to complete the
proof it remains to note that $T$ and $T'$ are conjugated in
$\SL(n,\R)$. \qed

\section{Closures of $D_I$-orbits when $\#I >2$}
\medskip

\subsection{Main Proposition} If $K$ is a $\mathrm{CM}$-field we denote by $F$
the totally real subfield of $K$ of index 2. In this case we denote
by $F_i$ the completion of $F$ with respect to the valuation $|\ . \
|_i$ on $K_i$  and by $\OO_F$ the ring of integers of $F$. We put
$A_F = \underset{i}{\prod}F_i$.

In this section $I = \{1, \cdots, l\}$ where $3 \leq l \leq r$.

\begin{prop}\label{prop4} Let
$h = (e, \cdots, \underbrace{u_l^-(\beta)u_l^+(\alpha)}_{l}, \cdots, e)
\in G$ where $u_l^-(\beta) = \left(
\begin{array}{cc} 1&0\\ \beta&1\\ \end{array}
\right)$, $u_l^+(\alpha) = \left(
\begin{array}{cc} 1&\alpha\\ 0&1\\ \end{array}
\right)$, $\alpha \in K^*$ and $\beta \in K_l$. The following
assertions hold:
\begin{enumerate}
\item[(a)] If $K$ is not a $\mathrm{CM}$-field then $\overline{D_I\pi(h)} =
G/\Gamma$;
\item[(b)] Let $K$ be a $\mathrm{CM}$-field and $d_{\alpha}$ be an element in $D$ such
that $d_{\alpha}^2 = \left(
\begin{array}{cc} \alpha&0\\ 0&\alpha^{-1}\\ \end{array}
\right)$. Then $\overline{D_{I,\R}\pi(h)} \supset
d_{\alpha}G_{\R}d_{\alpha}^{-1}\pi(e)$ and
$d_{\alpha}G_{\R}d_{\alpha}^{-1}\pi(e)$ is closed.
\end{enumerate}
\end{prop}

In order to prove the proposition we need the following lemma.

\begin{lem}\label{lem3} Let $K$ be a $\mathrm{CM}$-field and $\alpha \in
K^*$. Then
$$
\overline{F_l\alpha + \OO} = A_F\alpha + \OO.
$$
\end{lem}

{\bf Proof.} Let $n$ be a positive integer such that $n \alpha \in
\OO$. By the classical strong approximation theorem $\overline{F_l +
\OO_F} = A_F$. Since $A_F \cap \OO = \OO_F$ we have that $A_F + \OO$
is closed in $A$ (cf.\cite[1.13]{Raghunathan}) and, therefore,
$$
\overline{F_l + \OO} = A_F + \OO.
$$
Put
$$
L = \overline{F_l\alpha + \OO n\alpha} = A_F\alpha + \OO n\alpha.
$$
Since $\OO n \alpha$ has finite index in $\OO$, $L \cap \OO$ is a
lattice in $L$. Hence $L + \OO$ is a closed subgroup of $A$ which,
in view of the definition of $L$, completes the proof. \qed

\textbf{Proof of Proposition \ref{prop4}.} Note that $U^+(A)\pi(e)$
is closed and homeomorphic to $A/\OO$. (We denote by $U^+(A)$ the
group of $A$-points of the upper unipotent subgroup of $G$.) This
implies that $u_l^+(K_l)\pi(e)$ is dense in $U^+(A)\pi(e)$ and, when
$K$ is a $\mathrm{CM}$-field, it follows from Lemma \ref{lem3} that
$u_l^+(F_l\alpha)\pi(e)$ is dense in the closed set
$U^+(A_F\alpha)\pi(e)$.

Further the proof proceeds in several steps.

\medskip

{\it Step 1.} As in the formulation of Proposition \ref{integers1},
let  $H$ be the closure of the projection of $\OO^*$ into $K_l^*
\times \cdots \times K_r^*$. Denote by $p_j: A^* \rightarrow K_j^*$,
$l \leq j \leq r$, the natural projections. We will consider the
case (a) (when $K$ is not a $\mathrm{CM}$-field) and the case (b)
(when $K$ is a $\mathrm{CM}$-field) in a parallel way. Using
Propositions \ref{integers1}(a) and \ref{integers2}, for every
positive integer $m$ we fix in $H^{\circ}$ a compact neighborhood
$H_m$ of $1$ with the following properties: (i) $1 - \frac{1}{m} <
|p_j(x)|_j < 1 + \frac{1}{m}$ for all $j \geq l$ and all $x \in H_m$
and, (ii) $p_l(H_m) = \{e^{(a_m + \imath b_m)t}: t \in
[-\frac{1}{m}, \frac{1}{m}]\}$, where $\imath = \sqrt{-1}$ and $a_m$
and $b_m$ are real numbers such that $b_m \neq 0$ (resp. $b_m = 0$
and $a_m \neq 0$) if $K_l = \C$ and we are in case (a) (resp. if
otherwise). In view of Proposition \ref{integers1}(b) there exists a
sequence $y_n \in \OO^*$ such that $y_n \in \OO_F^*$ in case (b),
$|p_l(y_n)|_l
> n$ and $1 - \frac{1}{n} < |p_j(y_n)|_j < 1 + \frac{1}{n}$ for all $j
> l$.

\medskip

{\it Step 2.} Denote
$$
L_{mn} = \{x^2 : x \in y_n H_m\}.
$$
Let $W_\varepsilon$ be the $\varepsilon$-neighborhood of 0 in $A$
and $W_{\varepsilon,F}$ be the $\varepsilon$-neighborhood of 0 in
$A_F$. We claim that given $m$ for every $\varepsilon > 0$ there
exists a constant $n_\circ$ such that if $n > n_\circ$ then
\begin{equation}
\label{formula-2} A = W_\varepsilon + p_l(L_{mn}) + \OO
\end{equation}
in case (a), and
\begin{equation}
\label{formula-1} A_F = W_{\varepsilon,F} + p_l(L_{mn}) + \OO_F
\end{equation}
in case (b).

Note that the projections of $K_l$ into $A/\OO$ and of $F_l$ into
$A_F/\OO_F$ are dense and equidistributed. Since $|p_l(y_n)|_l > n$
this implies the claim in case (b) and in case (a) when $K_l = \R$.

Consider the case (a) when $K_l = \C$. If $\theta \in [0,2\pi)$ we
put $\R_{\theta} = e^{\imath\theta}\R$ and if $a < b$ we put
$[a,b]_{\theta} = e^{\imath\phi}[a,b]$ where $\R$ stands for the
subfield of reals in $K_l$. Since $\overline{K_l + \OO} = A$ it is
easy to see that for almost all $\theta \in [0, 2\pi)$ we have that
$\overline{\R_\theta + \OO} = A$ and, moreover, given $\varepsilon
> 0$ there exists $c_\varepsilon > 0$ such that if $b-a >
c_\varepsilon$ then
$$
A = W_\varepsilon + z + [a,b]_{\theta} + \OO, \ \forall z \in A.
$$
Now let $p_l(y_n) = r_n e^{\imath\frac{\psi_n}{2}}$ and $\psi_n
\underset{n}{\rightarrow} \psi$. Since the real $b_m$ in the
definition of $H_m$ is different from 0 there exists
$\frac{\theta}{2} \in (-\frac{1}{m}, \frac{1}{m})$ such that
$\overline{\R_{\theta + \psi} + \OO}= A$. Remark that since $r_n
\rightarrow +\infty$ the curvatures at the points of the plane curve
$p_l(L_{mn}) \subset \C$ are tending uniformly to 0 when $n
\rightarrow \infty$. Therefore for every real $\beta > 0$ end every
$\varepsilon > 0$ there exist a positive integer $n_\circ$ such that
for every $n > n_\circ$ there exists a $z \in K_l$ such that the
points of the segment $z + [0, \beta]_{\theta + \psi}$ are
$\varepsilon$-close to $p_l(L_{nm})$. This implies the claim.

\medskip

{\it Step 3.} Since $d(\xi^{-1})\pi(e) = \pi(e)$ for every $\xi \in
\OO^*$ we have that $(e, \cdots,
u_l^-(\xi^{-2}\beta)u_l^+(\xi^{2}\alpha),d_{l+1}(\xi)^{-1},\cdots,d_r(\xi)^{-1})\pi(e)$
belongs to $D_I\pi(h)$ (respectively, $D_{I,\R}\pi(h)$) if $K$ is
not (respectively, is) a $\mathrm{CM}$-field. Put
$$ X_{mn} \bydefn \{(e, \cdots,
\underbrace{u_l^-(x^{-2}\beta)u_l^+(x^{2}\alpha)}_{l},
\cdots,d_r(x)^{-1})\pi(e): x \in y_nH_m\}
$$
Since $y_nH_m \cap \OO^*$ is dense in $y_nH_m$,
\begin{equation}
\label{formula2}
X_{mn} \subset \overline{D_I\pi(h)}
\end{equation}
in case (a), and
\begin{equation}
\label{formula3} X_{mn} \subset \overline{D_{I,\R}\pi(h)}
\end{equation}
in case (b). Using the commensurability of $\OO$ and $\OO\alpha$ we
deduce from (\ref{formula-2}) and (\ref{formula-1}) that for every
$m$
\begin{equation}
\label{formula4} \overline{\bigcup_n p_l(L_{mn}\alpha) + \OO} = A
\end{equation}
in case (a) and
\begin{equation}
\label{formula5} \overline{\bigcup_n p_l(L_{mn}\alpha) + \OO} = A_F
\alpha +\OO
\end{equation}
in case (b). On the other hand, it follows from the definitions of
$H_m$ and $y_n$ that for every $\delta > 0$ there exists $c_\delta$
such that if $\min\{m, n\} > c_\delta$ then $|x^{-2}\beta|_l <
\delta$ and $||x|_j - 1| < \delta$ for all $x \in y_nH_m$. Now it
follows from (\ref{formula2}), (\ref{formula3}), (\ref{formula4}) and (\ref{formula5})
that $U^+(A) \pi(e) \subset \overline{D_I\pi(h)}$ in case (a) and
$U^+(A_F\alpha)\pi(e) \subset \overline{D_{I,\R}\pi(h)}$ in case
(b).

\medskip

{\it Step 4.} Let $B_1^+$ and $B_{1,\R}^+$ be the upper triangular
subgroup of $G_1$ and $G_{1,\R}$, respectively. In view of {\it Step
3}, $B_1^+\pi(e) \subset \overline{D_I\pi(h)}$ in case (a) and
$d_{\alpha}B_{1,\R}^+d_{\alpha}^{-1}\pi(e) \subset \overline{D_{I,\R}\pi(h)}$ in case (b).
Note that $B_1^+$ and $d_{\alpha}B_{1,\R}^+d_{\alpha}^{-1}$ are epimorphic subgroups of $G_1$
and $d_{\alpha}G_{1,\R}d_{\alpha}^{-1}$, respectively. It follows from
\cite[Theorem 1]{Shah-Weiss} that $\overline{B_1^+\pi(e)} =
\overline{G_1\pi(e)}$ and $\overline{d_{\alpha}B_{1,\R}^+d_{\alpha}^{-1}\pi(e)} =
\overline{d_{\alpha}G_{1,\R}d_{\alpha}^{-1}\pi(e)}$. Suppose we are in case (b).
It is easy to see that  $d_{\alpha}^{-1}\Gamma d_{\alpha}$
contains a congruence subgroup of $\Gamma$. Therefore
$d_{\alpha}^{-1}\Gamma d_{\alpha}$ and $\Gamma$
are commensurable and since $G_{\R}\pi(e)$ is closed
$d_{\alpha}G_{\R}d_{\alpha}^{-1}\pi(e)$ is too. Using, for example,
Borel's density theorem \cite{Raghunathan} one sees that
$\overline{G_1\pi(e)} = G/\Gamma$ and
$\overline{d_{\alpha}G_{1,\R}d_{\alpha}^{-1}\pi(e)} =
d_{\alpha}G_{\R}d_{\alpha}^{-1}\pi(e)$. Therefore
$\overline{D_I\pi(h)} = G/\Gamma$ in case (a) and
$\overline{D_{I,\R}\pi(h)} \supset
d_{\alpha}G_{\R}d_{\alpha}^{-1}\pi(e)$ in case (b). \qed

\subsection{Proofs of Theorem \ref{thm2} and Corollary \ref{cor2}.}
It is enough to prove Theorem \ref{thm2} for $I = \{1,2,3\}$. We may
(and will) assume that $g_i \in G_{i,K}, \ i \in I$. By the theorem
hypothesis either $g_1g_2^{-1} \notin \mathcal{N}_{G_K}(D_K)$ or
$g_2g_3^{-1} \notin \mathcal{N}_{G_K}(D_K)$ (see Proposition
\ref{lem}). Suppose that $g_1g_2^{-1} \notin
\mathcal{N}_{G_K}(D_K)$. In view of Proposition \ref{prop3} there
exists an element $\pi(g') \in \overline{D_I\pi(g)}, g' = (g_1',
\cdots, g_r')$, such that $g_i' \in G_K$ if $1 \leq i \leq 3$, $g_1'
g_2'^{-1} \in \mathcal{N}_{G_K}(D_K)$, $g_1'g_3'^{-1} \notin
\mathcal{N}_{G_K}(D_K)$ and $g_i' = g_i$ if $i > 3$. Clearly, if $n
\in {D_I}$ and $k \in G_K$ then ${D_I\pi(g')}$ is dense if and only
if ${D_I\pi(ng'k)}$ is dense (see Lemma \ref{lem2}(c)). Therefore we
may assume without loss of generality that $\overline{D_I\pi(g)}$
contains an element $\pi(h)$ where $h$ is as in the formulation of
Proposition \ref{prop4}. Now Theorem \ref{thm2} follows from
Proposition \ref{prop4}(a).

Let us prove Corollary \ref{cor2}. By Moore's ergodicity theorem
\cite{Z}, $D_{I,\R}$ is ergodic on $G/\Gamma$. Therefore there
exists an $y \in G/\Gamma$ such that $D_{I,\R}y$ is dense in
$G/\Gamma$. By Theorem \ref{thm2}, $\overline{D_I\pi(g)} =
G/\Gamma$. Therefore there exists a compact $M \subset D_I$ such
that $M\overline{D_{I,\R}\pi(g)} = G/\Gamma$. Let $y = mz$, where $m
\in M$ and $z \in \overline{D_{I,\R}\pi(g)}$. Then
$$
\overline{D_{I,\R}\pi(g)} \supset m^{-1}\overline{D_{I,\R}y} =
G/\Gamma
$$
which completes the proof. \qed

\medskip

\subsection{Proof of Theorem \ref{thm3}.} Recall that $I = \{1,
\cdots, l\}$, $l \geq 3$. Choose $g = (e,\cdots,
\underbrace{u_l^{+}(\alpha)u_l^{-}(\beta)}_{l}, \cdots, e)$ where
$\alpha \in K \setminus F$, and $\beta \in F^*$. We will prove that
$x = \pi(g)$ is the point we need. First, remark that
$u_l^{+}(\alpha)u_l^{-}(\beta) = tu_l^{-}(\beta_1)u_l^{+}(\alpha_1)$
where $t \in D_{l,K}$, $\beta_1 \in K$ and $\alpha_1 =
\frac{\alpha}{1+\alpha\beta}$. Hence $\alpha_1 \in K \setminus F$.
Let $d_{\alpha_1} \in D$ be such that $d_{\alpha_1}^2 = \left(
\begin{array}{cc} \alpha_1&0\\ 0&\alpha_1^{-1}\\ \end{array}
\right)$. Applying twice Proposition \ref{prop4}(b) we obtain that
\begin{equation}
\label{formula6} \overline{D_{I,\R}x} \supset G_{\R}\pi(e) \bigcup
d_{\alpha_1}G_{\R}d_{\alpha_1}^{-1}\pi(e).
\end{equation}
Note that the orbits $G_{\R}\pi(e)$ and
$d_{\alpha_1}G_{\R}d_{\alpha_1}^{-1}\pi(e)$ are closed and proper.

Since $G_{\R}\pi(e) \supset U^-(A_F)\pi(e)$ and
$d_{\alpha_1}G_{\R}d_{\alpha_1}^{-1}\pi(e) \supset
U^+(A_F\alpha_1)\pi(e)$ we have that
$$
\overline{D_{I,\R}x} \subset \{u_l^+(\mu\alpha)G_{\R}\pi(e): {0\leq
\mu\leq1}\} \bigcup
\{tu_l^-(\nu\beta_1)d_{\alpha_1}G_{\R}d_{\alpha_1}^{-1}\pi(e):
{0\leq \nu\leq1}\},
$$
where $\mu$ and $\nu \in F_l$. This implies (i).

Let us prove (ii). Using Proposition \ref{integers1} we can choose a
sequence $\xi_i \in \OO_F^*$ such that for every $j \geq l$ the
projection of $\xi_i$ into $F_j$ converges to some $x_j \in F_j^*$
and $x_l$ is not an algebraic number. Put
\begin{eqnarray*}
y = (e, \cdots, u_l^{+}(x_l^{2}\alpha)u_l^{-}(x_l^{-2}\beta),
d_{l+1}(x_{l+1}^{-1}), \cdots, d_r(x_{r}^{-1}))\pi(e).
\end{eqnarray*}
Then
$$
y = \underset{i}{\lim} \ d_I(\xi_i)x \in \overline{D_{I,\R}x}.
$$
Let us show that $y \notin D_{I,\R}x$. Otherwise, there exist
elements $d \in D_l$ and $m \in G_{l,K}$ such that
$du_l^{+}(x_l^{2}\alpha)u_l^{-}(x_l^{-2}\beta) =
u_l^{+}(\alpha)u_l^{-}(\beta)m$. Since
$u_l^{+}(\alpha)u_l^{-}(\beta)m \in G_{l,K}$ the lower right
coefficient of $du_l^{+}(x_l^{2}\alpha)u_l^{-}(x_l^{-2}\beta)$
belongs to $K$. This implies that $d \in D_{l,K}$ and that
$x_l^{2}\alpha \in K$ which contradicts our choice of $x_l$, proving
the claim.

Let $H$ be a subgroup of $G$ such that $H \supset D_{I,\R}$ and $Hy$
be closed. It is easy to see that
$$
x = \underset{i}{\lim} \ d_I(\xi_i^{-1})y.
$$
In view of (\ref{formula6}), $H$ contains both  $G_{\R}$ and
$d_{\alpha_1}G_{\R}d_{\alpha_1}^{-1}$. Since $\alpha_1 \in K
\setminus F$, we obtain $A = A_F + A_F\alpha_1$ and $A = A_F + A_F\alpha_1^{-1}$. Therefore, $H
\supset U^+(A) \cup U^-(A)$. Hence $H = G$ which proves (ii).

In order to prove (iii), suppose on the contrary that
$$
\overline{D_{I,\R}x} \setminus D_{I,\R}x \subset
\underset{i}{\bigcup}H_ix_i
$$
where $\{H_i\}$ is a countable family of proper subgroups of $G$ and
$H_ix_i$ are closed orbits. Then each $H_i$ is closed. Let $y$ be as
in the formulation of (ii). It follows from the Baire category
theorem that there exists $H_j$ such that $D_{I,\R} \subset H_j$ and
$y \in H_jx_j$. But the latter contradicts (ii). \qed
\subsection{Proof of Corollary \ref{cor3}.} The fact that $\mathcal{N}_{G}(D_I)G_K$
is not equal to $\bigcap_{i \in I}\big(\mathcal{N}_{G}(D_i)G_K\big)$
is easy to prove. In view of Theorem \ref{q+}(b) the orbit $D_I\pi(g)$ is
locally divergent if and only if $g \in \bigcap_{i \in I}(\mathcal{N}_{G}(D_i)G_K)$ and in view of Proposition \ref{lem} if
$g \in \mathcal{N}_G(D_I)G_K$ then $\overline{D_I\pi(g)}$ is an
orbit of a torus. Suppose that  $g \in \bigcap_{i \in I}\big(\mathcal{N}_{G}(D_i)G_K\big)\setminus \mathcal{N}_{G}(D_I)G_K$.
Let $g = (g_1, \cdots, g_r)$. There exist $i$ and $j \in I$, $i \neq j$,  such that $g_ig_j^{-1}$ does not normalise the diagonal
subgroup of $\SL(2)$.

We have seen in \S4.1 when $\#I = 2$ and in \S5.2
and \S5.3 when $\#I > 2$ that in this case $\overline{D_I\pi(g)}$ is
not an orbit of a torus. \qed

\section{A Number-theoretic Application}
\label{spit forms}

{\it In this section we prove Theorem \ref{application}.} We use the
notation preceding the formulation of the theorem.

We identify the elements from $G/\Gamma$ with the lattices in $A^2$
obtained via the injective map $g\Gamma \mapsto g\OO^2$. This map is
continuous and proper with respect to the quotient topology on
$G/\Gamma$ and the topology of Chabauty on the space of lattices in
$A^2$.

The group $G_K$ is acting on $K[X,Y]$ by the law
$$(\sigma p)(X,Y) =
p(\sigma^{-1}(X,Y)), \forall \sigma \in G_K, \ \forall p \in K[X,Y],
$$
where $\sigma^{-1}(X,Y) = (m_{11}X + m_{12}Y, m_{21}X + m_{22}Y)$,
$\sigma^{-1} = \left(\begin{array}{cc} m_{11}&m_{12}\\ m_{21}&m_{22}\\ \end{array} \right)$.
By the theorem hypothesis $f_i(X,Y) = l_{i,1}(X,Y) \cdot
l_{i,2}(X,Y)$ where $l_{i,1}$ and $l_{i,2} \in K[X,Y]$ are linearly
independent over $K$ linear forms. There exist $g_i \in G_{i,K}$
such that $f_i(X,Y) = \alpha_i(g_i^{-1}f_0)(X,Y)$ where $\alpha_i
\in K^*$ and $f_0$ is the form $X\cdot Y$. We may (and will) suppose
that $\alpha_i = 1$ for all $i$. Since the forms $f_i, \ 1 \leq i
\leq r$ are not proportional, $g = (g_1, \cdots, g_r)$ does not
belong to $\mathcal{N}_G(D)G_K$. Therefore $D\pi(g)$ is a locally
divergent non-closed orbit (Theorem \ref{q+}(b)).

Let $r > 2$. Fix $a = (a_1, \cdots, a_r) \in A$ and choose $h \in G$
such that $he_1 = (a,1)$ where $e_1$ is the first vector of the
canonical basis of the free $A$-module $A^2$. According to Theorem
\ref{thm2}, $D\pi(g)$ is a dense orbit. Therefore there exist $d_n
\in D$ and $\gamma_n \in \Gamma$ such that $\underset{n}{\lim}\
d_ng\gamma_n = h$. Put $z_n = \gamma_n e_1$. Then $z_n \in \OO^2$
and
$$
\underset{n}{\lim}\ f(z_n) = \underset{n}{\lim}\ f_0(d_n g\gamma_n
e_1) = f_0(\underset{n}{\lim}\ (d_ng\gamma_n(e_1))) = f_0(a,1) = a,
$$
which proves the part (a) of the theorem.

Let $r = 2$. We will prove the inclusion
\begin{equation}
\label{inclusion} \overline{f(\OO^2)} \subset
\underset{j=1}{\overset{s}{\bigcup}} \phi^{(j)}(\OO^2) \bigcup (K_1'
\times \{0\}) \bigcup (\{0\} \times K_2') \bigcup f(\OO^2),
\end{equation}
where $\phi^{(j)}$, $K_1'$ and $K_2'$ are as in the formulation of
the theorem. Let $a = (a_1, a_2) \in \overline{f(\OO^2)} \setminus
f(\OO^2)$. There exists a sequence $z_n = (\alpha_n,\beta_n)$ in
$\OO^2$ such that $a = \underset{n}{\lim}\ f(z_n)$ and $f(z_n) \neq
0$ for all $n$. Let $a_1 \neq 0$. (The case $a_2 \neq 0$ is
analogous.) If $f_2(z_n) = 0$ for infinitely many $n$ then it is
easy to see that $a \in K_1' \times \{0\}$. From now on we suppose
that $f_2(z_n) \neq 0$ for all $n$.

Let $g = (g_1,g_2) \in G$ be such that $g_i(X,Y) = (l_{i1}(X,Y),
l_{i2}(X,Y))$, $i \in \{1,2\}$. We choose sequences $s_n \in K_1^*$
and $t_n \in K_2^*$ such that
\begin{eqnarray}
\label{applic1} \left\{
\begin{array}{ccc}
\underset{n}{\lim} \ s_nl_{11}(z_n)&=&a_{11}\\
\underset{n}{\lim} \ s_n^{-1}l_{12}(z_n)&=&a_{12}
\end{array}
\right. \ \text{and} \ \ \ \left\{
\begin{array}{ccc}
\underset{n}{\lim} \ t_nl_{21}(z_n)&=&a_{21}\\
\underset{n}{\lim} \ t_n^{-1}l_{22}(z_n)&=&a_{22}
\end{array}
\right.
\end{eqnarray}
where $a_{11}, a_{12} \in K_1$, $a_{21}, a_{22} \in K_2$, $a_1 =
a_{11}\cdot a_{12}$ and $a_2 = a_{21}\cdot a_{22}$.

If $a_2 = 0$ we choose $t_n$ in such a way that
\begin{equation}
\label{lim-1} a_{21} = a_{22} = 0.
\end{equation}
We have
\begin{equation}
\label{lim0} \underset{n}{\lim} \ d(s_n,t_n)g(z_n) =
(\mathbf{a_1},\mathbf{a_2})
\end{equation}
where $\mathbf{a_1} = (a_{11}, a_{12}) \in K_1^2$ and $\mathbf{a_2}
= (a_{21}, a_{22}) \in K_2^2$.

Shifting $g$ from the left by an element from
$\mathcal{N}_{G_K}(D_{K})$ if necessary, we reduce the proof to the
case when $|s_n|_1 \rightarrow \infty$ and $|t_n|_2 \leq 1$. There
exist $\mu$ and $\nu \in K$ such that
$$
l_{22} = \mu l_{11} + \nu l_{12}.
$$
We have
\begin{multline*}
0 < |\mathrm{N}_{K/\Q}(l_{22}(z_n))| = |l_{22}(z_n)|_1 \cdot
|l_{22}(z_n)|_2 = \\ = {|s_n|_1}\cdot{|t_n|_2}\cdot|\mu
s_n^{-1}l_{11}(z_n) + \nu s_n^{-1}l_{12}(z_n)|_1\cdot
|t_n^{-1}l_{22}(z_n)|_2.
\end{multline*}
Since $\{\mathrm{N}_{K/\Q}(l_{22}(z_n))\}$ is a discrete subset of
$\R$ which does not contain $0$, in view of (\ref{applic1}), we
obtain that
\begin{equation}
\label{lim1}
\underset{n}{\liminf} {|s_n|_1}\cdot {|t_n|_2}
> 0
\end{equation}
and that
$$
|a_{22}|_2 = \underset{n}{\lim} |t_n^{-1}l_{22}(z_n)|_2 \neq 0.
$$
The latter contradicts (\ref{lim-1}). Hence $a_2 \neq 0$.

Let us prove that
\begin{equation}
\label{Bruhat} g_1g_2^{-1} \in B^-_KB^+_K.
\end{equation}
First we need to show that
\begin{equation}
\label{lim2} \underset{n}{\limsup} {|s_n|_1}\cdot {|t_n|_2} <
\infty.
\end{equation}
There exist $\mu'$ and $\nu' \in K$ such that
$$
l_{11} = \mu' l_{21} + \nu' l_{22}.
$$
Then
\begin{eqnarray*}
& & 0 < |\mathrm{N}_{K/\Q}(l_{11}(z_n))| = |l_{11}(z_n)|_1 \cdot
|l_{11}(z_n)|_2 = \\
& & = |s_n|_1^{-1}\cdot|t_n|_2^{-1}\cdot|s_nl_{11}(z_n)|_1 \cdot
|\mu' t_nl_{21}(z_n) + \nu' t_nl_{22}(z_n)|_2.
\end{eqnarray*}
Now (\ref{lim2}) follows from the inequality $|t_n|_2 \leq 1$ and
(\ref{applic1}).

Suppose on the contrary that $g_1g_2^{-1} \notin B^-_KB^+_K$.
Therefore $g_1g_2^{-1} \in \omega B^+_K$. Shifting $g$ from the left
by a suitable element from $D_K$ we reduce the proof to the case
when $g_1g_2^{-1} = \omega u$, where $u = \left(
\begin{array}{cc} 1&\alpha\\ 0&1\\ \end{array} \right)$.
In view of (\ref{lim1}), (\ref{lim2}), Lemma \ref{lem0} and Lemma
\ref{lem2} we can find a sequence $\xi_n \in \OO^*$ and a converging
to $a \in A^*$ sequence $a_n \in A^*$ such that $(s_n,t_n) =
\xi_na_n$ and $d(\xi_n)g_2\OO^2 = g_2\OO^2$. Using (\ref{lim0}) we
see that $d(\xi_n)g(z_n)$ converges to some
$(\mathbf{b_1},\mathbf{b_2}) \in A^2$ where $\mathbf{b_1} = (b_{11},
b_{12}) \in K_1^2$ and $\mathbf{b_2} = (b_{21}, b_{22}) \in K_2^2$.
(Recall that $A^2$ is identify with $K_1^2 \times K_2^2$.) An easy
computation shows that
$$
d(\xi_n)g(z_n) = (h_n, e)\mathbf{w_n}
$$
where $h_n = \left(
\begin{array}{cc} 0&\xi_n^2\\ -\xi_n^{-2}&-\alpha\\ \end{array}
\right)$ and $\mathbf{w_n} = d(\xi_n)g_2(z_n) = (\beta_n, \gamma_n)
\in g_2\OO^2$. So, $\big((\xi_n^2 \gamma_n, -\xi_n^{-2}\beta_n -
\alpha \gamma_n), (\beta_n, \gamma_n)\big) \rightarrow
(\mathbf{b_1},\mathbf{b_2})$ which implies that $(\xi_n^2 \gamma_n,
\gamma_n)$ converges to $(b_{11},b_{22})$ in $A$. But
$$|\xi_n^2
\gamma_n|_1\cdot |\gamma_n|_2 = |\xi_n^2|_1\cdot
|\mathrm{N}_{K/\Q}(\gamma_n)|.$$ Hence
$$
\underset{n}{\lim}\ |\xi_n^2|_1\cdot |\mathrm{N}_{K/\Q}(\gamma_n)| =
|b_{11}|_1\cdot |b_{22}|_2
$$
which is a contradiction because $|\xi_n^2|_1 \rightarrow \infty$
and $\underset{n}{\liminf}\ |\mathrm{N}_{K/\Q}(\gamma_n)| > 0$. This
completes the prove of (\ref{Bruhat}).

In view of Proposition \ref{prop2}(b), there exists a subsequence of
$d(s_n,t_n)\pi(g)$ converging to an element from
$\underset{j=1}{\overset{s}{\cup}}D\pi(h_j), \ 2 \leq s \leq 4$
where $h_j \in \mathcal{N}_{G_K}(D_{K})$  (see Corollary
\ref{cor1}). So, there exists $d \in D$ such that
$(\mathbf{a}_1,\mathbf{a}_{2}) \in dh_j\OO^2$, $1 \leq j \leq s$.
Hence $a \in \underset{j=1}{\overset{s}{\cup}}\phi^{(j)}(\OO^2)$
where $\phi^{(j)} = h_j^{-1}f_0$. It is clear that the quadratic forms
 $\phi^{(j)}$ are $2 \times 2$ nonproportional.
This completes the proof of
(\ref{inclusion}).

The inclusion inverse to (\ref{inclusion}) is easy to prove. Let $c
= \phi^{(j)}(z)$ where $z \in \OO^2$. We have $h_j =
\underset{n}{\lim} \ t_n g \sigma_n$ for some $t_n \in D$ and
$\sigma_n \in \Gamma$. Therefore
$$
\phi^{(j)}(z) = \underset{n}{\lim} \ f_0(t_n g \sigma_n(z)) =
\underset{n}{\lim} \ f(\sigma_n(z)) \in \overline{f(\OO^2)}.
$$
It remains to prove that $(K_1' \times \{0\}) \bigcup (\{0\} \times
K_2') \subset \overline{f(\OO^2)}$. It is enough to prove that if
$(x,y) \in K_1^2$ and $f_2(x,y) = 0$ then $(f_1(x,y),0) \in
\overline{f(\OO^2)}$. Suppose that $l_{21}(x,y) = 0$. Since $l_{11}$
and $l_{12}$ are linear combinations of $l_{21}$ and $l_{22}$ we get
that $f_1(x,y) = c\cdot l_{22}(x,y)^2$ where $c$ is a constant. Note
that the projection of the set $\{l_{22}(z): z \in \OO^2, l_{21}(z)
= 0\}$ into $K_1$ is dense. Therefore $(f_1(x,y),0) \in
\overline{f(\OO^2)}$. By similar reasons if $l_{22}(x,y) = 0$ then
$f_1(x,y) = d\cdot l_{21}(x,y)^2 \in \overline{f(\OO^2)}$, where $d$
is a constant. Note that $K_1' = c\{\alpha^2: \alpha \in K_1\} \cup
d\{\alpha^2: \alpha \in K_1\}$ and that, since $f_1$ and $f_2$ are
not proportional, $c$ and $d$ can not be simultaneously equal to
zero. This readily implies that $K_i' = \C$ if $K_i = \C$, and that
$K_i' = \R, \R_{-}$ or $\R_{+}$ if $K_i = \R$. The proof is
complete.\qed

\section{Concluding remarks}
\label{conclusion}

The elements $h_i$ in the formulation of Theorem
\ref{thm1} can be explicitly described in terms of $g$. This becomes
clear from the proof of this theorem in \S 4.1. Here we will give an
example of an orbit $D_I\pi(g), \ I = \{1,2\},$ such that the
boundary of its closure consists of {\it four} different closed
orbits.

In the next proposition we suppose that $G/\Gamma$ is a Hilbert
modular space of rank $r = 2$. Let $K_v$ be the completion of $K$
with respect to a non-archimedean valuation $v$ of $K$. Since $K$ is
dense in $K_v$ we may (and will) choose $\alpha$ and $\beta \in K$
such that $\alpha \cdot \beta \neq 0$, $\alpha \cdot \beta \neq 1$,
$|\alpha|_v > 1$ and $|\beta|_v < 1$.

\begin{prop}
\label{s=4} With the above notation and assumptions, let $g =
(g_1,g_2) \in G$ where $g_1 = \left(
\begin{array}{cc} 1&0\\ \alpha &1\\ \end{array}
\right)$ and $g_2 = \left(
\begin{array}{cc} 1&\beta\\ 0&1\\ \end{array}
\right)$. Then
$$\overline{D\pi(g)} \setminus {D \pi(g)} =
\underset{i=1}{\overset{4}{\cup}}D\pi(h_i),$$ where $D\pi(h_i)$ are
pairwise different, closed, non-compact orbits.
\end{prop}

{\bf Proof.} Since the coefficients of the matrix $g_1 g_2^{-1}$ are
different from 0, all pairs $(\sigma_1, \sigma_2) \in \{0,1\}^2$ are
admissibles and, in view of (\ref{limit2}), we need to prove that
the closed orbits $D(\omega^{\sigma_1}, \omega^{\sigma_2})
\pi(h_{\sigma_1,\sigma_2})$ are pairwise different. We have seen in
the course of the proof of Theorem \ref{thm1} that
$D(\omega^{\sigma_1}, \omega^{\sigma_2}) \pi(h_{\sigma_1,\sigma_2})
\neq D(\omega^{\sigma_1'}, \omega^{\sigma_2'})
\pi(h_{\sigma_1',\sigma_2'})$ if $(\sigma_1, \sigma_2) = (0,0)$ or
$(1,1)$ and $(\sigma_1', \sigma_2') = (0,1)$ or $(1,0)$. It remains
to show that $D\pi(h_{0,0}) \neq D\pi(\omega h_{1,1})$ and
$D(\omega, 1) \pi(h_{1,0}) \neq D(1,\omega)\pi(h_{0,1})$.

Using (\ref{boundary}) we see that $h_{0,0} = e$ and modulo
multiplication from the left by an element from $D_K$, $\omega
h_{1,1}$ is equal to $\left(
\begin{array}{cc} \frac{1}{1-\alpha\beta}&\frac{\beta}{1-\alpha\beta}\\ \alpha &1\\ \end{array}
\right)$. Since $\alpha \notin \OO$ we conclude that $D\pi(h_{0,0})
\neq D\pi(\omega h_{1,1})$.

Modulo multiplication from the left by an element from $D_K$,
$h_{1,0}$ (respectively, $h_{0,1}$) is equal to $\left(
\begin{array}{cc} 1&\frac{1}{\alpha}\\ 0&1\\ \end{array}
\right)$ (respectively, $\left(
\begin{array}{cc} 1&0\\ \frac{1}{\beta}&1\\ \end{array}
\right)$). If $D(\omega, 1) \pi(h_{1,0}) = D(1,\omega)\pi(h_{0,1})$
then
$$
\frac{\xi^2 \beta + \alpha}{\alpha\beta} \in \OO
$$
for some $\xi \in \OO^*$. This leads to contradiction because, in
view of the choice of $\alpha$ and $\beta$,
$$
\frac{|\xi^2 \beta + \alpha|_v}{|\alpha\beta|_v} =
\frac{1}{|\beta|_v} > 1.
$$
Therefore the boundary of $D_I\pi(g)$ consists of four pairwise
different closed orbits. \qed

\medskip

Remark that most of the results of this paper remain valid with very
small changes in the $S$-adic setting, that is, when $G$ is a
product of $\SL(2,K_v)$, where $K_v$ is the completion of a number
field $K$ with respect to a place $v$ belonging to a finite set $S$
of places of $K$ containing the archimedean ones. For instance, the
proofs of the analogs of Theorems \ref{thm1} and
\ref{application}(b) remain valid in this context with virtually no
changes. When $K$ is not a $\mathrm{CM}$-field, the analog of
Theorem \ref{thm2} remains true with very small modifications if $K
= \Q$ or if $I$ contains an archimedean place. For instance, Theorem
\ref{thm2} remains true for action of maximal tori (that is, when $D
= D_I$). The analog of Theorem \ref{thm2} in the general case (for
arbitrary $K$ and $I$) is more delicate and will be treated later.
Also, one can find tori orbits with non-homogeneous closures for
many spaces $G/\Gamma$ with $G \neq \SL_n$. This will be treated
elsewhere too.

\medskip

\emph{Acknowledgements:} I would like to thank both Yves Benoist and
Nimish Shah for the profitable discussions. I am grateful to Elon
Lindenstrauss for a useful discussion on Conjecture B.

I wish to thank Grisha Margulis for his useful remarks on a
preliminary version of the paper.

I am grateful to the organizers of the Oberwolfach Workshop
"Homogeneous Dynamics and Number Theory", July 4th-July 10th 2010,
for giving me the opportunity to report the results of this paper.

\end{document}